\documentclass[11pt, a4]{amsart}

\usepackage[english]{babel}
\usepackage{amsmath,enumerate, amsfonts, amssymb,amsthm, xypic}
\usepackage[dvips]{graphics} 

\newtheorem{thm}{Theorem}[section]
\newtheorem{lemma}[thm]{Lemma}
\newtheorem{prop}[thm]{Proposition}
\newtheorem{cor}[thm]{Corollary}

\newenvironment{demo}{\noindent{\it Proof.}\,}{\begin{flushright}
    \,$\Box$ \smallskip \end{flushright}}

\theoremstyle{definition}
\newtheorem{defi}[thm]{Definition}

\newtheorem{ex}[thm]{Example}
\newtheorem{rmk}[thm]{Remark}

\newcommand{\eps}{\varepsilon}

\DeclareMathOperator{\ord}{ord}
\DeclareMathOperator{\jac}{jac}
\DeclareMathOperator{\mult}{mult}
\DeclareMathOperator{\Reg}{Reg}
\DeclareMathOperator{\Sing}{Sing}
\DeclareMathOperator{\inte}{int}

\title{ Motivic invariants of Arc-Symmetric sets and Blow-NASH
  Equivalence}

\author{Goulwen Fichou}

\address {D\'epartement de Math\'ematiques,
Universit\'e d'Angers, 2 bd Lavoisier, 49045 Angers Cedex, France}

\email{fichou@tonton.univ-angers.fr}
\begin{document}

\begin{abstract} We define invariants of the blow-$Nash$ equivalence of
  real analytic function germs, in a similar way that the motivic zeta
  functions of Denef \& Loeser \cite{DL1}. As a key ingredient, we
  extend the virtual Betti numbers, which were known for real
  algebraic sets \cite{MCP}, as a generalized Euler characteristics for
  projective constructible
  arc-symmetrics sets. Actually we prove more: the virtual Betti numbers
  are not only algebraic invariant, but also $Nash$-invariant of
  arc-symmetric sets. Our zeta functions enable to sketch the blow-$Nash$
  equivalence classes of Brieskorn polynomials of two variables.
\end{abstract}

\maketitle

\section*{Introduction}
In the study of real analytic function germs, the choice of a good
equivalence relation is an important question. Whereas a topological equivalence is too coarse and a
$C^1$-equivalence too fine, blow-analytic equivalence, a notion
introduced by T.-C. Kuo in 1985 ( see \cite{Kuo}, and \cite{FKK}
for a survey) seems to behave better, especially with respect to
finiteness properties. In this paper we will focus on a particular case of
blow-analytic equivalence, say blow-$Nash$ equivalence. Let
$f,g:(\mathbb R^d,0) \longrightarrow (\mathbb R,0)$ be real analytic
function germs; $f$ and $g$ are said to be blow-$Nash$ equivalent if there
exist two algebraic modifications $\pi_f:\big( M_f,\pi_f^{-1}(0) \big) \longrightarrow (\mathbb R^d,0)$
and $\pi_g:\big( M_g,\pi_g^{-1}(0)\big) \longrightarrow (\mathbb R^d,0)$, and a
$Nash$-isomorphism, that is an analytic isomorphism with
semi-algebraic graph, $\phi :\big( M_f,\pi_f^{-1}(0)\big) \longrightarrow
\big( M_g,\pi_g^{-1}(0)\big) $ which respects the multiplicity of the
jacobian determinants
of $\pi_f$ and $\pi_g$ and which
induces a homeomorphism $h:(\mathbb R^d,0) \longrightarrow (\mathbb
R^d,0)$ such that $f=g\circ h$. Here by a modification of $f$ we mean
a proper birational map which is an isomorphism
over the complement of the zero locus of $f$.
One
can define such an equivalence relation on germs of real analytic
varieties, and S. Koike proved (\cite{Koi1,Koi2}) that finiteness
properties hold in this case. However in the case of germs of
functions, the question of moduli is still an open question.

A common issue for blow-analytic equivalence and blow-$Nash$
equivalence is to prove that, when it is the case, two given function
germs are not equivalent. The difficulty rises in the lack of
invariants known for these equivalence relations. Up to now, just one
kind of invariants have been known:
the Fukui invariants. To an analytic function germ
$f$, the Fukui invariants associate the set of possible orders $n$ of series $f \circ
\gamma(t)= a_nt^n+\cdots, a_n\neq 0$, for $\gamma: (\mathbb R,0) \longrightarrow (\mathbb R^d,0)$
an analytic arc (\cite{Fukui,IKK}). There exists also a version of the
Fukui invariants related to the sign of $f$.

Using motivic integration conbined with the construction of a
computable motivic invariant for arc-symmetrics sets, the virtual Betti
numbers, we introduce in this paper zeta
functions $Z(T), Z^{\pm}(T)$ of a real analytic function germ that
belong to $\mathbb Z[u,u^{-1}][[T]]$, and take into account
not only the orders of series $f \circ
\gamma(t)$ but also the geometry of the sets $\chi_n(f)$ of arcs $\gamma$ that
realize a given order $n$ (for
precise definitions, see section \ref{defzeta}). These zeta functions
are similar to the motivic zeta functions of Denef \& Loeser \cite{DL1}. We prove that
actually
they are invariants for blow-$Nash$ equivalence. The proof is directly
inspired by the work of Denef \& Loeser via their formulae for the
zeta functions in terms of a modification of the zero locus of the
analytic function germ (propositions \ref{DLnaive}, \ref{DLmono}). It
uses the powerfull machinery of motivic integration, a theory introduced by
M. Kontsevitch in 1995 \cite{Kont} and developped by J. Denef
and F. Loeser \cite{DL1,DL2,DL3,DL4}, in  particular the fundamental change of variable formula
(\ref{chgvar}).

In order to dispose of computable invariants, motivic integration
requires computable measures, or in other words Generalized Euler
Characteristics. A Generalized Euler
Characteristic is an additive and multiplicative invariant defined on
the level of the Grothendieck group of varieties. In our setting of
the blow-$Nash$ equivalence, we need
invariants of the Zariski constructible sets over real algebraic varieties
$\chi_n(f)$ ( real algebraic variety is the sense of \cite{BCR} ), and we ask it to be respected by $Nash$-isomorphisms. It leads
naturally to the category of $Nash$-varieties, and more generally of
arc-symmetric sets.

Arc-symmetric sets have been introduced in 1988 by K. Kurdyka \cite{KK1} in
order to study ``rigid components'' of semi-algebraic
sets. Arc-symmetric sets enabled him to prove Borel theorem
\cite{KK2} which says that injective endomorphisms of real algebraic
sets are surjective. With a slightly different definition of arc-symmetric sets
A. Parusi\'nski has proved the same result by using the fact that
these sets form a constructible category \cite{AP}.

In section \ref{sect2}, we give conditions on an invariant defined
on connected components of compact nonsingular real algebraic
varieties such that it extends to a additive invariant on the constructible category of
arc-symmetric sets. Additive means that
$\chi(A)=\chi(B)+\chi(A \setminus B)$ for $B \subset A$ an inclusion
of arc-symmetric sets. Let us stress the fact that the unique such
additive invariant known
up to now in the real case is the classical Euler characteristic with compact support,
and as a matter of fact it is the unique Generalized Euler
Characteristic for semi-algebraic sets up to homeomorphism
\cite{Q}. As a fundamental example, we prove that the Betti numbers with $\mathbb Z
_2$-coefficient defined on connected components $A$ of compact
nonsingular real
algebraic sets by $b_k(A)= \dim H_k(A,\mathbb Z _2)$, give such an
invariant $\beta_k$ on arc-symmetric sets for each
$k \in \mathbb N$ (corollary \ref{cor-main}), called $k$-Virtual Betti
number. We make them
multiplivative by putting $\beta(A)=\sum_{k= 0}^{\dim A} \beta_k(A)u^k
\in \mathbb Z[u]$, called virtual Poincar\'e polynomial of $A$. This invariant is 
different from the classical Euler characteristic with compact support, and in particular it enables
to distinguish homeomorphic real algebraic sets which are not
isomorphic. Moreover it respects dimension as put in light by the
formula $\deg \big(\beta(A)\big)=\dim(A)$ (see remark
\ref{degdim}) whereas the Euler characteristic with compact support
identifies the dimension.  These numbers have been proven to be additive
invariant of real algebraic varieties recently by C. McCrory and Adam
Parusi\'nski in \cite{MCP}; in this paper we extend the virtual Betti
numbers to the more general context of arc-symmetric sets, and we
prove the invariance not only under algebraic isomorphisms but also under
$Nash$-isomorphisms (see \ref{isoNash}).

Note that the virtual Betti numbers over real algebraic sets have been
introduced independly by C. McCrory
and A. Parusi\'nski \cite{KP}, and by B. Totaro \cite{Tot}. Moreover S. Koike
and A. Parusi\'nski have defined in the same way zeta functions, by
using the classical Euler characteristic with compact support, which are invariant
for blow-analytic equivalence (\cite{KP}). The advantage of our zeta
functions, whose invariance is proven only for blow-$Nash$
equivalence, is that the Virtual Betti numbers have a better behaviour
with respect to algebraicity and analyticity that the classical Euler
characteristic with compact support which is merely topological.

In the first section we study the behaviour of arc-symmetric sets with
respect to closure, irreductibility, nonsingularity and overall
resolution of singularities (proposition \ref{ecl}). In the second one we show that, under
certain conditions, an invariant defined over connected components
of compact nonsingular algebraic sets can be extended to an additive
invariant on all arc-symmetric sets (cf theorem \ref{main}). Here by
invariant we mean up to algebraic isomorphisms (see \ref{defisoalg}). This
result requires the use of the Weak Factorization Theorem
(\cite{weak,wlo}) and of the desingularisation Theorem of H. Hironaka \cite{HIRO}, and is in the spirit of the description of the naive
Grothendieck group of $k$-varieties, with $k$ a field of
characteristic zero, given by F. Bittner in \cite{FB}. In
particular this part leads to the existence of the Generalized Euler
Characteristic $\beta$, defined from the Betti numbers, for arc-symmetric sets.
Then we extend the invariance of $\beta$ in section \ref{sec3}, not
only to algebraic isomorphisms of arc-symmetric sets but also to
$Nash$-isomorphisms (for a precise definition see \ref{isoNash}).
In the final part we define zeta functions of a real analytic function
germ and we
prove, by the way of motivic integration, the main theorem of this part: $Z(T), Z^{\pm}(T)$ are invariants
for blow-$Nash$ equivalence. In particular we state Denef \&
Loeser formulae for the zeta functions (proposition \ref{DLnaive},\ref{DLmono}). Finally we apply these
invariants to sketch the blow-$Nash$
equivalence classes of Brieskorn polynomials in two variables, and to give some examples in
three variables.

{\bf Acknowledgement.} I am greatly indebted to Adam Parusi\'nski, my thesis
advisor, for his help during this work at the university of Angers.

\section{Arc-Symmetric sets}
\subsection{Arc-Symmetric sets and closure}

We fix a compactification of $\mathbb R^n$, for instance $\mathbb R^n
\subset \mathbb P^n$.
\begin{defi} Let $A \subset \mathbb P^n$ be a semi-algebraic set. We say
  that $A$ is arc-symmetric if one of the two equivalent conditions
  holds:
\begin{enumerate}
\item for every real analytic arc $\gamma
  :]-\epsilon,\epsilon[ \longrightarrow  \mathbb P^n$ such that $\gamma
  (]-\epsilon,0[) \subset A$ there exists $\epsilon ' > 0$ such that $\gamma
  (]0,\epsilon '[) \subset A$,
\item for every real analytic arc $\gamma
  :]-\epsilon,\epsilon[ \longrightarrow  \mathbb P^n$ such that $\inte
  \gamma ^{-1}(A) \neq \varnothing$ there exists $x_0, \ldots, x_n \in
  ]-\epsilon,\epsilon[ $ such that $\gamma (]-\epsilon,\epsilon[
  \setminus \{ x_0, \ldots, x_n \}) \subset A$.
\end{enumerate}
\end{defi}

Remark that an arc-symmetric set need not to be an analytic variety
(cf \cite{KK1}, example 1.2).

This definition is the one of A. Parusi\'nski \cite{AP}. Note that a
closed arc-symmetric set is necessarily compact. This definition
  differs from the one of K. Kurdyka \cite{KK1} who only considers closed
  arc-symmetric sets in $\mathbb R^n$. One can think our arc-symmetric
  sets as projective constructible arc-symmetric sets.
Remark that the arc-symmetric sets form a constructible
category $\mathcal {AS}$ of semi-algebraic sets:
\begin{itemize}
\item $\mathcal {AS}$ contains the algebraic sets,
\item $\mathcal {AS}$ is stable under set-theoretic operations $\cup,\,
  \cap,\, \setminus$,
\item $\mathcal {AS}$ is stable by inverse images of $\mathcal
  {AS}$-map (i-e whose graph is in $\mathcal {AS}$)  and images of injective $\mathcal {AS}$-map,
\item each $A \in \mathcal {AS}$ has a well-defined fundamental class
  with coefficients in $\mathbb Z _2$.
\end{itemize}

In particular there is a notion of closure in $\mathcal {AS}$ (we
refer to \cite{AP} for a proof).

\begin{figure}
\caption{The Whitney umbrella.}\label{f1}
$$\resizebox{4cm}{!}{\rotatebox{0}{\includegraphics{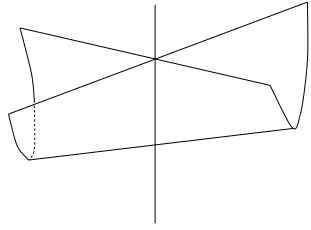}}}$$
\end{figure}

\begin{prop} Every $A \in \mathcal {AS}$ admits a smallest
  arc-symmetric set, denoted by $\overline A ^{\mathcal {AS}} $,
  containing $A$ and closed in $\mathbb P^n$. 
\end{prop}
\begin{rmk} Even if $A \in \mathcal {AS}$, we do not have that $\overline A \in
  \mathcal {AS}$. Consider for $A$ the regular part of the Whitney umbrella $zx^2=y^2$ (see figure
  \ref{f1}). The closure of $A$ in $\mathcal {AS}$ is the entire Whitney
  umbrella.
\end{rmk}

One can consider irreducible arc-symmetric sets: $A \in \mathcal {AS}$ is
irreducible if $A=B \cup C$, with $B$ and $C$ closed in $A$ and
arc-symmetric, implies that $ B \subset C$ or $C \subset B$. Remark
that an irreducible arc-symmetric set is not necessarily connected
with our definition of arc-symmetric sets (as an example consider a
hyperbola in the plane). Nethertheless, as proven in \cite{KK1}, an
arc-symmetric set $A$ admits a unique decomposition as a finite union of
irreducible arc-symmetric sets closed in $A$. Note that $\mathcal
{AS}$-closure has a good behaviour with respect to irreducibility:

\begin{prop}\label{irr} If $A \in \mathcal {AS}$ is irreducible, then so is $\overline A
 ^{\mathcal {AS}}$.
\end{prop} 
\begin{demo} Suppose $\overline A
 ^{\mathcal {AS}}=B \cup C$ with $B$ and $C$ arc-symmetric and closed
 in $\overline A ^{\mathcal {AS}}$. Then $B$ and $C$ are closed and
 $A$ splits in $A=A \cap  A ^{\mathcal {AS}}= (A \cap B) \cup (A \cap
 C)$, with $A \cap B$ and $A \cap
 C$ arc-symmetric and closed in $A$. But $A$ is irreducible so, for
 instance, $A \cap B \subset A \cap C$. Then $A=A \cap C$, so $A
 \subset C$ and $\overline A ^{\mathcal {AS}} = C$ because $C$
 is arc-symmetric and closed. 
\end{demo}

\begin{rmk} The hypothesis ``$B$ and $C$ closed'' is
  essential. Let $A$ be the regular part of the Whitney umbrella. The $\mathcal {AS}$-closure
  of $A$ is the entire Whitney umbrella which is the disjoint union of two
  arc-symmetric sets: $A$ and the vertical line.
\end{rmk} 
 Recall that the dimension of an arc-symmetric set is its dimension as
 a semi-algebraic set, and also the dimension of its Zariski
 closure in the projective space \cite{KK1} (here by real algebraic
 variety we mean is the sense of \cite{BCR}). In particular, if $A \in \mathcal {AS}$ one has $\dim A=\dim \overline A
^{\mathcal {AS}}= \dim \overline A ^Z$.

\begin{prop}\label{dimen} Let $A \in \mathcal {AS}$. Then $\overline A ^{\mathcal {AS}}=A \cup  \overline {\overline A
   \setminus A }^{\mathcal {AS}}$. In particular $\dim \overline A
 ^{\mathcal {AS} } \setminus A
   < \dim A$.
\end{prop} 
\begin{demo} As a union of arc-symmetric sets $F=A \cup  \overline {\overline A
   \setminus A }^{\mathcal {AS}}$ is arc-symmetric. Moreover $F =A \cup (\overline A
  \setminus A)  \cup (\overline {\overline A
   \setminus A }^{\mathcal {AS}})= \overline A  \cup \overline {\overline A
   \setminus A }^{\mathcal {AS}}$ thus $F$ is closed. So we have
 proved that $\overline A ^{\mathcal {AS}} \subset F$.
Moreover $\overline A \subset \overline A ^{\mathcal {AS}}$ because
$\overline A ^{\mathcal {AS}}$ is closed, thus $\overline {\overline A
   \setminus A }^{\mathcal {AS}} \subset \overline A ^{\mathcal
   {AS}}$, hence $F \subset \overline A ^{\mathcal
   {AS}}$. Consequently $F=\overline A ^{\mathcal
   {AS}}$.
\end{demo}

We can adapt proposition 1.5 of \cite{KK2} to our definition of
arc-symmetric sets.

\begin{prop}\label{sub-irr} Let $A \in \mathcal{AS}$ be irreducible, and
  $B \subset A$ be a closed arc-symmetric subset of $A$ of the same
  dimension. Then $B=A$.
\end{prop}
%
%

\begin{demo} $A$ can be decomposed into the union of two arc-symmetric
  sets closed in $A$: $A=B \cup (\overline {A \setminus
    B} ^{\mathcal {AS}} \cap A)$. Then, by irreducibility of $A$,
  either $B \subset \overline {A \setminus B} ^{\mathcal {AS}} \cap A$
  or $\overline {A \setminus B} ^{\mathcal {AS}} \cap A \subset B$. In
  the second case $B=A$, and in the first one $B \subset (\overline {A
    \setminus B} ^{\mathcal {AS}} ) \setminus (A \setminus B)$. But
  this can not happen because
  the dimension of this arc-symmetric set is stricly less than $\dim
  B$ by proposition \ref{dimen}.
\end{demo}


\subsection{Nonsingular arc-symmetric sets}
 Let us define a nonsingular arc-symmetric set with relation
 to its Zariski closure in the projective space.
\begin{defi} An arc-symmetric set $A$ is nonsingular if $A \cap
  \textrm{Sing} (\overline A ^Z)=\varnothing$.
\end{defi}

\begin{lemma}\label{lis-conn} A nonsingular and connected arc-symmetric set is irreducible.
\end{lemma} 

\begin{demo} Let $C \in \mathcal {AS}$ be nonsingular and
  connected, and $A \subset C$ a closed connected arc-symmetric set of the same
  dimension. We are going to prove that $A=C$.

Denote by $A_0$ the semi-algebraic set consisting of the part of
maximal dimension of $A$, and put $A_1=\overline {A_0}$. Then $A_1$
is a closed semi-algebraic subset of $C$, contained in $A$. But $A_1$ is also open in
$C$: take $a \in A_1$; there exists an open ball $D$ with center $a$
included in $C$. Then $\dim D \cap A_0=\dim A$, therefore $D \subset
A$ because $A$ is arc-symmetric (one can fill $D$ with analytic arcs
whose interiors intersect $A_0$ in non empty sets). Moreover $D
\subset A_0$, by definition of $A_0$, and then $A_1$ is an open
neighbourhood of $a$. Finally $A_1 $ is a connected component of $C$, so
$A_1=C$. But $A_1 \subset A$, therefore $A=C$.
\end{demo}

Let us state the definition of an isomorphism between arc-symmetric sets.
\begin{defi}\label{defisoalg} Let $A,B \in \mathcal {AS}$. Then $A$ is isomorphic to
  $B$ if and only if there exist Zariski open subsets $U$ and $V$ in $ \overline A
  ^{\mathcal {Z}}$ and $\overline B ^{\mathcal {Z}}$ containing $A$
  and $B$, and an algebraic
  isomorphism $\phi: U \longrightarrow V$ such that $\phi(A)=B$.
\end{defi}


The following proposition says that closed and nonsingular
arc-symmetric sets are very similar to compact nonsingular real algebraic
varieties.

\begin{prop}\label{comp-lis} Let $A \in  \mathcal {AS}$ be compact and
  nonsingular. Then $A$ is isomorphic to a union of connected
  components of some compact nonsingular real algebraic variety.
\end{prop}
\begin{demo} Let $X=\overline A ^{\mathcal {Z}}$ be the Zariski closure
  of $A$ in the projective space, and $\pi : \widetilde X \longrightarrow X$ a resolution of
  singularities of $X$. Remark that $\dim A=\dim X=\dim \widetilde
  X$ and that $A$ is isomorphic to the subset
  $\pi^{-1}(A)=\widetilde A$ of $\widetilde X$ because $A \subset \Reg(X)$ and
  $\Reg(X)$ is a Zariski open subset of $X$ isomorphic to $ \pi ^{-1}\big(
  \Reg(X)\big) \subset  \widetilde X$. Denote by
  $\widetilde X=\bigcup _{i \in I} C_i$ the decomposition of
  $\widetilde X$ in connected components. Each $C_i, i \in I$ is a
  closed and nonsingular arc-symmetric set, hence irreducible by
  proposition \ref{lis-conn}. Therefore $\widetilde A \cap  C_i=C_i$
  or $\varnothing$ because of proposition \ref{sub-irr}, and
  $\widetilde A$ is a union of connected components of $\widetilde X$
  as claimed.
\end{demo}


\subsection{Arc-Symmetric sets and resolution of singularities}
The following proposition is just an adaptation of Theorem 2.6
of \cite{KK1} to our definition of arc-symmetric sets.
\begin{prop}\label{ecl} Let $A \in \mathcal {AS}$ be irreducible. Let $X$ be a
  real algebraic variety containing $A$ with $\dim X=\dim A$, and $\pi :
  \widetilde X \longrightarrow X$ a resolution of singularities for
  $X$ (cf \cite{HIRO}). There exists a unique connected component $ \widetilde A$ of $
  \widetilde X$ such that $ \pi (\widetilde A)=\overline {\Reg(A)}$.
\end{prop}
%
%
%

\begin{demo}  Let $\widetilde {A_0}$ be an irreducible arc-symmetric component of
  dimension $\dim A$ of $\pi ^{-1}(A)$ ( such an $\widetilde {A_0}$ exists because
  $\dim \pi ^{-1}(A)=\dim A$). Then $\widetilde {A_0}$ is contained in some
  connected component $\widetilde A$ of $\widetilde X$. Actually $\widetilde {A_0}
  \subset \overline {\widetilde {A_0}} ^{\mathcal
   {AS} }$, which is irreducible because so is $\widetilde {A_0}$
 (proposition \ref{irr}) and  which is closed, is connected (proposition\ref{comp-lis}). Now $\overline {\widetilde A_0} ^{\mathcal
   {AS} }$ is included in some connected component of $\widetilde X$, and is equal to it by
  proposition \ref{sub-irr}. We can put  $\widetilde A=\overline {\widetilde {A_0}} ^{\mathcal
   {AS} }$.

 Let us show that $\pi (\widetilde
  A)=\overline{ \Reg(A)}$. In fact, it suffices to show that $\dim \overline{\pi (\widetilde
  A)}^\mathcal {AS} \setminus \pi (\widetilde A) < \dim A$, what will
be done in the next lemma. Indeed, on
one hand $\overline{\pi (\widetilde
  A)}^\mathcal {AS} =\overline A ^{\mathcal
   {AS} }$ by proposition \ref{sub-irr}, so $\dim \overline A ^{\mathcal
   {AS} } \setminus \pi (\widetilde A) < \dim A$. Now $\Reg(A) \cap \big(\overline A ^{\mathcal
   {AS} } \setminus \pi (\widetilde A)\big)$ is an open subset of $\overline A ^{\mathcal
   {AS} }$ of dimension stricly less than $\dim A$, so $\Reg(A) \cap \overline A ^{\mathcal
   {AS} } \setminus \pi (\widetilde A)=\varnothing$. This implies $\Reg(A)
\subset \pi (\widetilde A)$. 
On the other hand, if $E$ denotes the exceptionnal
divisor of the resolution, one has $\pi (\widetilde A \setminus E)
\subset \Reg(\overline A ^{\mathcal
   {AS} })$. However $\Reg(\overline A ^{\mathcal
   {AS} }) \subset  \overline {\Reg(A)}$ because $\dim \overline A ^{\mathcal
   {AS} } \setminus A < \dim A$ by proposition \ref{dimen}. 

Thus we have the following inclusions $\pi (\widetilde A \setminus E)
\subset  \overline {\Reg(A)}  \subset  \pi (\widetilde A)$ which gives
the conclusion by taking the closure.
\end{demo}

\begin{lemma} Let $A$ and $\widetilde A$ be as in the proof of
  proposition \ref{ecl}. Then $$\dim \overline{\pi (\widetilde A)}^\mathcal {AS} \setminus \pi (\widetilde A) < \dim A.$$
\end{lemma}
\begin{demo} Let us show
that $ \overline{\pi (\widetilde A)}^\mathcal {AS} \subset F= \pi
(\widetilde A) \cup \overline{\pi (E)}^{\mathcal {AS}}$. Remark that if $F$ is closed and arc-symmetric the result is proved. As $\pi$ is proper $\pi
(\widetilde A)$ is closed and so is $F$. Now, let $\gamma
  : ]-\epsilon,\epsilon[ \longrightarrow  \mathbb P^n$ be an real
  analytic arc such that $\inte \gamma ^{-1}(F) \neq \varnothing$. Then
  either $\inte \gamma ^{-1}(\overline{\pi (E)}^{\mathcal {AS}} ) \neq \varnothing$ and $\gamma
  (]-\epsilon,\epsilon[) \subset \overline{\pi (E)}^{\mathcal {AS}} $,
  or $\inte \gamma ^{-1}(\pi (\widetilde A) \setminus \overline{\pi
    (E)}^{\mathcal {AS}} ) \neq \varnothing$. In the latter case there exists a unique analytic
  arc $\widetilde {\gamma} : ]-\epsilon,\epsilon[ \longrightarrow
  \mathbb P^n$ such that $\pi \circ \widetilde {\gamma}=\gamma$. One
  has $\inte \widetilde {\gamma}^{-1}(\widetilde A) \neq \varnothing$,
  so $\widetilde {\gamma} (]-\epsilon,\epsilon[) \subset \widetilde A
  $ anf finally $ \gamma (]-\epsilon,\epsilon[) \subset \pi (
  \widetilde A) \subset F$. Thus $F$ is arc-symmetric.
\end{demo}

\begin{rmk}\label{iso} Denote by $D$ the singular locus of $X$ and by $E$ the
  exceptionnal divisor of the resolution. Then $\pi : \widetilde A
  \setminus E \longrightarrow \overline A ^{\mathcal {AS}} \setminus
  D$ is an isomorphism of arc-symmetric
  sets (restriction of an algebraic ismorphism).

 If we add the assumption that $A$ is nonsingular, then the
  conclusion becomes $ \pi (\widetilde A)=\overline {A}$. Moreover
  $\pi : \widetilde A \setminus \big(\widetilde A \setminus \pi^{-1}(A)\big) \longrightarrow A$
 is an isomorphism of arc-symmetric sets, and $\dim \widetilde A
 \setminus \pi^{-1}(A) < \dim A$.
\end{rmk}
One can describe with precision the differences between $A$ and $ \pi
(\widetilde A)$. Actually the symmetric difference of $A$ and $ \pi
(\widetilde A)$ consists of a semi-algebraic set of dimension stricly
less than $\dim A$, and more precisely one has the following
proposition.
\begin{prop} Let $A$ and $\widetilde A$ be as in proposition \ref{ecl}. Then
$$A \setminus ( \pi (\widetilde A) \cap A)=\{ x\in \Sing(A),\,\, \dim_x A <
\dim A \}$$ and
$$ \pi (\widetilde A) \setminus  \big(A \cap \pi (\widetilde A)\big)=
\Reg(\overline A ^{\mathcal {AS} }) \setminus A \cup \{ x \in \Sing(
\overline A ^{\mathcal {AS} } \setminus A),\,\, \dim_x \overline A ^{\mathcal
   {AS} }= \dim A \}$$
$$= \{ x \in
\overline A ^{\mathcal {AS} } \setminus A,\,\, \dim_x \overline A ^{\mathcal
   {AS} }= \dim A \}.$$
\end{prop}
\begin{demo} For the first equality, remark that $A \setminus ( \pi
  (\widetilde A) \cap A) \subset \textrm{Sing}(A)$ and that
  $\textrm{Sing}(A)$ splits in $\{ x\in \Sing(A),\,\, \dim_x A <
\dim A \} \cup \{ x\in \Sing(A),\,\, \dim_x A =
\dim A \}$. But recall that $ \pi (\widetilde A)=\overline {\Reg(A)}$,
and so $\Sing(A) \cap \pi (\widetilde A)= \Sing(A) \cap \overline
{\Reg(A)}=\{ x\in \Sing(A),\,\, \dim_x A = \dim A \}$.

In the same way, $\pi (\widetilde A) \setminus  \big(A \cap \pi
(\widetilde A)\big) \subset \overline A ^{\mathcal {AS} } \setminus A$ and
$\overline A ^{\mathcal {AS} } \setminus A = \Reg(\overline A
^{\mathcal {AS} } \setminus A) \cup S_1 \cup S_2$ with
$S_1= \{ x\in \Sing \big( \overline A
^{\mathcal {AS} } \setminus A)\big),\,\, \dim_x \overline A ^{\mathcal
  {AS}} = \dim A\}$ and 
$S_2= \{ x\in \Sing \big(\overline A
^{\mathcal {AS} } \setminus A)\big),\,\, \dim_x \overline A ^{\mathcal
  {AS}} < \dim A\}$.
One has $S_2 \cap \overline {\Reg(A)}=\varnothing$ and $S_1 \subset \overline
{\Reg(A)}$. Moreover the remark \ref{iso} induces that $\Reg(\overline A
^{\mathcal {AS} } \setminus A) \subset \overline A ^{\mathcal {AS} }
\setminus D  \subset \pi (\widetilde A)$, hence the conclusion. 
\end{demo}

\begin{figure}

\caption{Resolution of the Whitney umbrella.}\label{f2}
$$\resizebox{8cm}{!}{\rotatebox{0}{\includegraphics{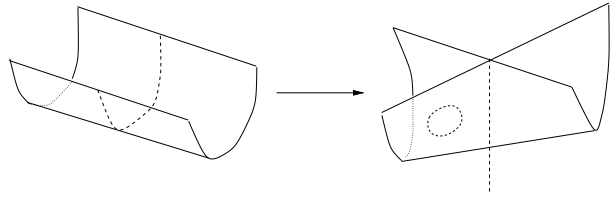}}}$$
\end{figure}
\begin{ex} Consider once more the Whitney umbrella $zx^2=y^2$. One can resolve the singularities by blowing-up along the
  $z$-axis. Let $\sigma (u,v,w)=(u,uv,w)$ denote the chart which
  contains the hole strict transform of the Whitney umbrella; then it
  has equation $w=v^2$, while the exceptional disisor has equation
  $u=0$ ( figure
  \ref{f2}).

On one
  hand, if $A$ is just the
  regular part minus the circle in dotted line, then $A \subset \pi (\widetilde
  A)$ and  $ \pi (\widetilde
  A) \setminus  \big(A \cap \pi (\widetilde A)\big)$ consists of the circle
  and the closed upper part of the vertical line. On the other and if $A$ is
  the entire Whitney umbrella, then $ \pi (\widetilde A ) \subset A$ and $A
  \setminus ( \pi (\widetilde A) \cap A)$ is the open bottom part of the
  vertical line.
\end{ex}

Let us finish this section by giving the well-known particular case of the ``blow-up'' of
a closed nonsingular arc-symmetric set. We state the proposition in
terms of connected component of nonsingular compact algebraic
varieties (recall proposition \ref{comp-lis}).

\begin{prop}\label{ecl-comp} Let $Y \subset X$ be compact nonsingular algebraic
  varieties such that $\dim Y < \dim X$, and $A \subset X$ be a connected
  component of $X$. Denote by $\pi : \widetilde X \longrightarrow X$
  the blow-up of $X$ along $Y$. Then $\pi$ is surjective and $\pi  ^{-1}(A)$ is a connected component of $\widetilde  X$.
\end{prop}

\section{Generalized Euler characteristic of arc-symmetric sets}\label{sect2}

\begin{defi} An additive map on $\mathcal {AS}$ with values
  in an abelian group is a map $\chi$ defined on $\mathcal {AS}$
  such that
\begin{enumerate}
\item for $A$ and $B$ arc-symmetric sets with $A$ isomorphic to $B$, then $\chi(A)=\chi(B)$,
\item for $B \subset A$ inclusion of arc-symmetric sets, then
  $\chi(A)=\chi(B)+\chi(A \setminus B)$.
\end{enumerate}
If moreover $\chi$ satisfies $\chi(A\times B)=\chi(A) \times \chi(B)$
for $A,B$ arc-symmetric sets, then we say that $\chi$ is a generalized Euler characteristic on $\mathcal {AS}$.
\end{defi}

\begin{rmk} The Euler characteristic with compact supports is a
  generalized Euler characteristic on $\mathcal {AS}$. Actually if we
  consider just semi-algebraic sets, with isomorphism replaced by
  homeomorphism, the Euler characteristic with compact support is the
  unique generalized Euler characteristic (see \cite{Q}).
However, for complex algebraic varieties, there exist a lot of such generalized
Euler characteristic, for example given by mixed Hodge structure (see
\cite{DL1, Loo}) 
\end{rmk}
The aim of this section is to give sufficiently good conditions on an
invariant $\chi$ over the closed (i-e compact) and nonsingular arc-symmetric
sets such that $\chi$ extends to an additive map on
$\mathcal {AS}$. We state the theorem in terms of connected
components of real algebraic varieties thanks to proposition
\ref{comp-lis}. The method is inspired from the one of
\cite{FB}, and the main result is:

\begin{thm}\label{main} Let $\chi$ be a map defined on connected
components of compact nonsingular real algebraic varieties with values in an abelian group and such that
\begin{description}
\item[P1] $\chi(\varnothing)=0$,
\item[P2] if $A$ and $B$ are connected components of compact nonsingular real algebraic
  varieties which are isomorphic as arc-symmetric sets, then $\chi(A)=\chi(B)$,
\item[P3] with notations and assumptions of proposition \ref{ecl-comp},
  $$\chi\big(\pi  ^{-1}(A)\big)-\chi\big(\pi  ^{-1}(A) \cap \pi ^{-1}(A \cap Y)\big)=\chi (A)- \chi (A \cap Y).$$

\end{description}

Then $\chi$ extends to an additive map defined on $\mathcal {AS}$.
\end{thm}

Before giving the proof of the theorem, let us give some consequences.
First this result enables us to give one example of such an additive
map by considering the homology with $\mathbb Z _2$ coefficients.
\begin{cor}\label{cor-main} Let $\beta _i$ be defined on a connected component $A$ of a
  compact nonsingular algebraic variety by $\beta _i (A)= \dim
  H_i(A,\mathbb Z _2)$, the Betti number mod-2 of $A$, for $i \in \mathbb
  N$. Then $\beta _i$ extends to an additive map defined on $\mathcal {AS}$.
\end{cor}
\begin{demo} Property P1 and P2 of theorem \ref{main} are clear. Let us prove property P3. With
  notations of P3, the exact sequences with coefficients in $\mathbb Z
  /2$ of the pairs $\big(\widetilde
  A,\widetilde A \cap \pi ^{-1}(B)\big)$ and $( A,B)$ give the following
  commutative diagram
 $$\xymatrix{\cdots  \ar[r]  &  H_{i-1}\big(\widetilde A,\widetilde A \cap \pi ^{-1}(B)\big)
   \ar[r] \ar[d] & H_{i}\big(\widetilde A \cap \pi ^{-1}(B) \big)  \ar[r]
   \ar[d]  & H_i (\widetilde A) \ar[r]\ar[d] &\cdots \\
       \cdots  \ar[r]&      H_{i-1}( A,B)   \ar[r]  & H_i( B)  \ar[r]  & H_i( A) \ar[r] &\cdots} $$
where the vertical arrows are induced by $\pi$. Note that $\pi _*
:H_{i-1}\big(\widetilde A,\widetilde A \cap \pi ^{-1}(B)\big) \longrightarrow
H_{i-1}( A,B)$ is an isomorphism because $\pi$ is a homeomorphism
between $ \widetilde A \setminus \widetilde A \cap \pi ^{-1}(B)$ and
$A \setminus B$, and that $\pi _* : H_i
(\widetilde A) \longrightarrow  H_i( A)$ is surjective because $\pi$
is of degre 1 and the varieties are $\mathbb Z  /2$-oriented. Now it
is an easy game to check that the following sequence
$$0 \longrightarrow H_{i}(\widetilde A \cap \pi ^{-1}(B) )
\longrightarrow  H_i( B) \oplus  H_i (\widetilde A) \longrightarrow
H_i( A) \longrightarrow 0$$ is exact, hence  
$$\beta_i(\widetilde A)-\beta_i\big(\widetilde A \cap \pi ^{-1}(B)\big)=\beta_i
(A)- \beta_i (B).$$
Then we can apply theorem \ref{main}.
\end{demo}

It turns out to be easy to
adapt theorem \ref{main} in order to obtain not only additive maps
but also generalized Euler characteristics.
\begin{thm}\label{main2} Let $\chi$ be as in theorem \ref{main}. Assume moreover
  that $\chi$ takes values in a commutative ring, and that for
  connected components of compact
  nonsingular real algebraic varieties $A,B$ the relation $\chi(A\times B)=\chi(A) \chi(B)$ holds. Then the extension of $\chi$ on
  $\mathcal {AS}$ of theorem \ref{main} is a generalized Euler characteristic.
\end{thm}

The following corollary is an immediate consequence of K\"unneth formula.
\begin{cor}\label{cor2} Let $\beta$ be defined on $A \in \mathcal
  {AS}$ by $$\beta(A)=\sum _{i= 0}^{\dim A}\beta _i(A)u^i.$$ Then
  $\beta$ is a generalized Euler characteristic on $\mathcal {AS}$.
\end{cor}

\begin{figure}
\caption{resolution of $C_1$ and $C_2$}\label{f4}
$$\resizebox{8cm}{!}{\rotatebox{0}{\includegraphics{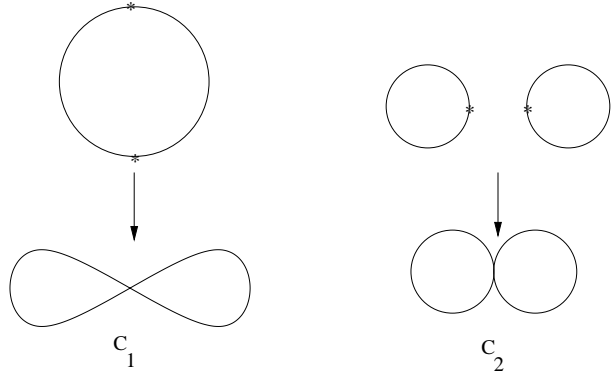}}}$$
\end{figure}

\begin{ex} 
\begin{enumerate}
\item If $\mathbb P^k$ denote the real projective space of dimension
  $k$, which is nonsingular and compact, then $\beta(\mathbb
  P^k)=1+u+\cdots+u^k$. Now, compactify $\mathbb R$ in $\mathbb P^1$
  by adding one point at the infinity. By additivity $\beta(\mathbb
  R)=\beta(\mathbb P^1)-\beta(point)=u$, and so $\beta(\mathbb R^k)=u^k$.
\item Let $W$ be the Whitney umbrella, and $L$ the line included in
  $W$. Therefore $\beta(W)=\beta(W \setminus L)+\beta(L)$ and
  $W\setminus L$ is isomorphic, by blowing-up, to the strict transform
  of $W$ minus a parabola $P$. So $$\beta(W \setminus L)=\beta(\mathbb
  A_{\mathbb R}^1 \times P)-\beta(P)=(\beta(\mathbb A_{\mathbb R}^1)-1)\beta(P)=(u-1)u.$$
Finally $\beta(W)=u^2$.
\item The real algebraic varieties of figure \ref{f4} (clearly, one
  can find reduced equations for such
  curves) are not isomorphic whereas they are homeomorphic. Indeed one can compute
  $\beta(C_1)$ and $\beta(C_2)$ by considering the resolutions given
  by blowing-up the singular point of $C_1$ and $C_2$. One finds
  $\beta(C_1)=u$ and $\beta(C_2)=2u-1$. Remark that the Euler
  characteristic (one can recover by evaluating $u$ at $-1$ in
  this example) does not distinguish these two curves.
\end{enumerate}
\end{ex}

\begin{rmk}\label{degdim} The invariant $\beta$ has the following
  property as we will see in the proof; let $A \in \mathcal
  {AS}$. Then
$$\dim(A)=\deg \big(\beta(A)\big).$$
\end{rmk}

Let us give the proof of theorem \ref{main} first, and later of theorem
\ref{main2}.

{\it Proof of theorem \ref{main}.} 
We prove theorem \ref{main} by induction on dimension; the rank $n$
inductive hypothesis claims that $\chi$ is defined on arc-symmetric
sets of dimension less than or equal to $n$, invariant under isomorphism of
arc-symmetric sets, and additive.

For $n=0$ the arc-symmetric sets are just finite unions of points and
the result is clearly true.
Assume that the inductive hypothesis is true for $n-1$. We prove the
result at rank $n$ in two steps:
\begin{enumerate}
\item if $\chi$ is an additive map on the nonsingular elements
  of $\mathcal {AS}$ of dimension less than or equal to $n$, then $\chi$ extends to an additive map on
  all arc-symmetric sets of dimension less than or equal to $n$,
\item if $\chi$ satisfies property P1,P2 and P3, then $\chi$ extends to
  an additive map on the nonsingular elements
  of $\mathcal {AS}$ of dimension less than or equal to $n$. 
\end{enumerate}
 
\textbf {Step 1} Let $A \in \mathcal {AS}$ of dimension $n$. There exists a
stratification $\overline A ^{\mathcal {Z}}=\bigcup _{S \in \mathcal
  S}S$ of $\overline A ^{\mathcal {Z}}$ with nonsingular algebraic
strata, that is the union is a disjoint union of locally closed
algebraic varieties (note in particular that we do not ask the strata
to be connected). Then $S \cap A$, for each $S \in \mathcal S$, is a nonsingular
arc-symmetric set, and $\chi(S \cap A)$ is defined. Put
$\chi(A)=\sum_{S \in \mathcal S}\chi(S \cap A)$. 
One has to check that $\chi(A)$ is well-defined and satisfies the
property of an additive function over $\mathcal {AS}$.

We show firstly that $\chi(A)=\sum_{S \in \mathcal S}\chi(S \cap A)$ in the
case where $A$ is nonsingular, by induction on the number of
elements in $\mathcal S$. Indeed take $N_0 \in \mathcal S$,
then $\chi(A)=\chi\big(A \setminus (A \cap N_0)\big) + \chi(A \cap N_0)$ and
$\chi\big(A \setminus (A \cap N_0)\big)=\sum_{S \in \mathcal S \setminus
  \{N_0\}}\chi(S \cap A)$ by induction, hence the result.

Now, if $\mathcal S_1$ and  $\mathcal S_2$ are two stratifications of
$\overline A ^{\mathcal {Z}}$, one can find a common refinement
$\mathcal S$. The independence in the nonsingular case induces that
$\sum_{S \in \mathcal S_1}\chi(S \cap A)=\sum_{S \in \mathcal S}\chi(S
\cap A)=\sum_{S \in \mathcal S_2}\chi(S \cap A)$, then $\chi$ does not
depend on the choice of the stratification.

Let us show that $\chi$ is additive; let $B,A \in \mathcal {AS}$ with $B
\subset A$ of dimension less than or equal to $n$. One can choose a stratification $\bigcup _{S \in \mathcal
  S}S$ of
$\overline A ^{\mathcal {Z}}$ such that $\overline B ^{\mathcal {Z}}$
and $\overline {A \setminus B} ^{\mathcal {Z}}$ are union of
strata. Then $$\sum_{S \in \mathcal S}\chi(S \cap B) + \sum_{S \in
  \mathcal S}\chi\big(S \cap (A \setminus B)\big)=\sum_{S \in \mathcal S}
\Big(\chi(S \cap B)+\chi\big(S \cap (A \setminus B)\big)\Big)$$
and $\chi\big(S \cap (A \setminus B)\big)+\chi(S \cap B)= \chi(S \cap A)$
because the strata are nonsingular, so
$\chi$ is additive.

\vskip 5mm

\textbf {Step 2}
The second step constitutes the heart of the work. Define $\chi$ over
the nonsingular arc-symmetric sets of dimension $n$ in the
following way:
\begin{description}
\item[D1] if $A=\bigcup _{i \in I} A_i$ denotes the decomposition
of $A$ into irreducible components, put $\chi(A)=\sum _{i \in I} \chi(A_i)$,
\item[D2] if $A \in \mathcal {AS}$ is nonsingular and irreducible, then
$\chi (A)=\chi (\widetilde A) -\chi \big(\widetilde A \setminus \pi^{-1}(A)\big)$, where
$\widetilde A$ is the connected component of a resolution of
singularities of $\overline A
^{\mathcal {Z} }$ given by remark \ref{iso} and $E$ the exceptionnal
divisor of the resolution.
\end{description}
 We have to prove that $\chi$ is well-defined, invariant under isomorphisms and additive
 over the nonsingular elements of $ \mathcal {AS}$.

The following lemma will be usefull in the sequel.
\begin{lemma}\label{indep} Let $A,B \in \mathcal{AS}$ be nonsingular, irreducible
  and isomorphic. Suppose that $ \overline A^{\mathcal Z}$
  and $ \overline B^{\mathcal Z}$ are nonsingular, and denote by $\widetilde A \subset \overline A^{\mathcal Z}$ and
  $\widetilde B \subset \overline B^{\mathcal Z}$ the
  connected components containing $A$ and $B$. Then $\chi (\widetilde A) -\chi (\widetilde A \setminus A)=\chi (\widetilde B) -\chi (\widetilde B \setminus B)$.
\end{lemma}  
\begin{demo} By
  definition of an isomorphism between arc-symmetric sets, we know
  that  $ \overline A^{\mathcal Z}$ and $ \overline B^{\mathcal Z}$
  are birationally equivalent, and the
  weak factorization theorem \cite{weak,wlo} factors this birational
  isomorphism in a succession of blow-ups and blow-downs. In
  particular we can assume that the birational isomorphism is just a
  blow-up $\pi:  \overline A^{\mathcal Z} 
  \longrightarrow  \overline B^{\mathcal Z} $ along a nonsingular
  variety $C$ such that $C \cap B=\varnothing$.
Note that $\pi ^{-1}(\widetilde B)=\widetilde A$ by proposition
  \ref{ecl-comp}.

Now $\chi(\widetilde B \setminus B)=\chi(\widetilde B \cap
C)+\chi\big(\widetilde B \setminus (B \cup C)\big)$ by the additivity inductive
hypothesis for $\dim \widetilde B \setminus B<\dim B$ by proposition \ref{dimen}. Moreover $\chi(\widetilde B)-\chi(\widetilde B \cap
C)=\chi(\widetilde A)-\chi\big(\widetilde A \cap \pi ^{-1}( C)\big)$ by
property P3, and $\chi\big(\widetilde B \setminus (B \cup
C)\big)=\chi\Big(\widetilde A \setminus \big(A \cup \pi ^{-1}( C)\big)\Big)$ by the
inductive hypothesis on invariance under isomorphism. Therefore
$\chi (\widetilde B) -\chi (\widetilde B \setminus B)=\chi(\widetilde A)-\chi\big(\widetilde A \cap \pi ^{-1}(
C)\big)-\chi\Big(\widetilde A \setminus \big(A \cup \pi ^{-1}( C)\big)\Big)$ which is
equal to $\chi (\widetilde A) -\chi (\widetilde A \setminus A)$ by the additivity inductive
hypothesis.
\end{demo}

\vskip 5mm
Let us check that the definition of $\chi$ for the nonsingular and
irreducible arc-symmetric sets of
dimension $n$ does not depend on the choice of the resolution of
singularities of remark \ref{iso}.

Let $A \in \mathcal {AS}$ be nonsingular and
  irreducible, and let $\pi _i : \widetilde {X_i} \longrightarrow \overline
  A^{\mathcal Z}$ for $i \in \{1,2\}$ be resolutions of singularities of $ \overline A^{\mathcal Z}$. Let $\widetilde
  {A_i}$ be the
  connected components of $\widetilde {X_i}$ given by proposition
  \ref{ecl}. One has to show that
$$\chi(\widetilde {A_1})-\chi\big(\widetilde {A_1} \setminus \pi_1^{-1}(A)\big)=\chi(\widetilde {A_2})-\chi\big(\widetilde {A_2} \setminus \pi_2^{-1}(A)\big).$$

But $\pi_1^{-1}(A)$ and $\pi_2^{-1}(A)$ are
  isomorphic irreducible nonsingular arc-symmetric sets because $\pi _i$
  is an isomorphism on a Zariski open subsets of $\widetilde {X_i}$
  containing $\pi_i^{-1}(A)$, for $i \in \{1,2\}$. Therefore lemma
  \ref{indep} applies and $\chi$ is well-defined.


\vskip 5mm
Now let us show that $\chi$ is invariant under isomorphisms of
arc-symmetric sets. The proof is very similar to the last one. 
 Let $A,B \in \mathcal{AS}$ be irreducible and isomorphic, which means that there exists Zariski open subsets $U$ and $V$ in $ \overline A
  ^{\mathcal {Z}}$ and $\overline B ^{\mathcal {Z}}$, and an algebraic
  isomorphism $\phi: U \longrightarrow V$ such that
  $\phi(A)=B$. Choose resolutions of singularities $\pi_A: \widetilde
  {X} \longrightarrow \overline A
  ^{\mathcal {Z}}$ and $\pi_B: \widetilde {Y} \longrightarrow \overline B ^{\mathcal {Z}}$ for $ \overline A
  ^{\mathcal {Z}}$ and $\overline B ^{\mathcal {Z}}$.
Then $ \pi_A^{-1}(A)$ and $\pi_B^{-1}(B)$ are isomorphic as
arc-symmetric sets, and then by lemma \ref{indep}
$\chi\big(\pi_A^{-1}(A)\big)=\chi\big(\pi_B^{-1}(B)\big)$. Moreover one has
$\chi\big(\pi_A^{-1}(A)\big)=\chi(A)$ because both are equal to $\chi(\widetilde A)-\chi\big(\widetilde A
  \setminus \pi_A^{-1}(A)\big)$, where $\widetilde A$ is the connected
  component of $\widetilde  {X}$ given by proposition \ref{ecl}. In the
  same way one has $\chi\big(\pi_B^{-1}(B)\big)=\chi(B)$, hence
  $$\chi(A)=\chi\big(\pi_A^{-1}(A)\big)=\chi\big(\pi_B^{-1}(B)\big)=\chi(B).$$

In the case where $A$ and $B$ are not irreducible, it suffices to
decompose $A$ and $B$ in irreducible components, and apply the property D1 because
an isomorphism between arc-symmetric sets respects the irreducible components.

\vskip 5mm

Finally, let us check that $\chi$ is additive. Let $B \subset A$ be an inclusion
of nonsingular arc-symmetric sets. Note that by definition of $\chi$
we need to prove the result only in the case where $A$ is irreducible.

 If $\dim
B=\dim A$, then $\overline B ^{\mathcal {AS}}=\overline A ^{\mathcal
  {AS}}$ by proposition \ref{sub-irr}, and $\overline B ^{\mathcal {Z}}=\overline A ^{\mathcal
  {Z}}$. Now choose a resolution of singularities $\pi :
\widetilde X \longrightarrow   \overline A ^{\mathcal {Z}}$ for
$\overline A ^{\mathcal {Z}}$. If $\widetilde A$ denotes the connected
component of $\widetilde X$ given by proposition \ref{ecl} for $A$,
then it is also the component associated to $B$, and therefore $\chi(B)=\chi(\widetilde A)-\chi\big(\widetilde A
  \setminus \pi^{-1}(B)\big)$. Now $\chi\big(\widetilde A  \setminus
  \pi^{-1}(B)\big)=\chi\big(\widetilde A \setminus
  \pi^{-1}(A)\big)+\chi\big(\pi^{-1}(A) \setminus \pi^{-1}(B)\big)$ by the
  inductive hypothesis on additivity. So
  $\chi(B)=\chi(A)-\chi(A\setminus B)$ because $\chi\big(\pi^{-1}(A)
  \setminus \pi^{-1}(B)\big)=\chi(A\setminus B)$ by the invariance under
  isomorphism in dimension smaller that $n$.

If $\dim B<\dim A$, choose a resolution of singularities $\pi :
\widetilde X \longrightarrow   \overline A ^{\mathcal {Z}}$ for
$\overline A ^{\mathcal {Z}}$. Then it is also a resolution of
singularities of $\overline {A\setminus B} ^{\mathcal {Z}}=\overline A ^{\mathcal {Z}}$. Then $\chi(A\setminus B)=\chi(\widetilde {A})-\chi\big(\widetilde {A} \setminus
\pi^{-1}(A\setminus B)\big)$. 
Now $$\chi \big(\widetilde {A} \setminus
\pi^{-1}(A\setminus B) \big)=\chi \Big(\big(\widetilde {A} \setminus \pi^{-1}(A)\big) \cup
\pi^{-1}(B) \Big)=\chi \big(\widetilde {A} \setminus
\pi^{-1}(A) \big)+\chi \big(\pi^{-1}(B) \big)$$ by the inductive
assumption, and once more
by the inductive assumption one has $\chi\big(\pi^{-1}(B)\big)=\chi(B)$. Finally $\chi(A\setminus B)=\chi(A)-\chi(B)$.

This achieves the proof of step 2, and thus the proof of theorem \ref{main}.
\begin{flushright}
$\Box$
\end{flushright}

As it was the case for the previous proof, we are going to prove
theorem \ref{main2} by induction on dimension. The following relations
will be usefull:
\begin{equation}\label{r1}
\chi(\sqcup _{i=1}^k A_i)=\sum _{i=1}^k \chi(A_i),
\end{equation}
where the union of the arc-symmetric sets $A_i, i=1,\ldots,k$ is disjoint, and 
\begin{equation}\label{r2}
\chi(A)=\chi(\widetilde A)-\chi(\widetilde A \setminus A),
\end{equation}
where $A$ is a nonsingular arc-symmetric set whose arc-symmetric
closure $\widetilde A$ is nonsingular.

{\it Proof of theorem \ref{main2}.} Put, as an inductive hypothesis at rank $n$, that $\chi$ is
  multiplicative for all arc-symmetric sets of dimension stricly less
  than or equal to $n$.

Remark that we can restrict our attention to the nonsingular case
because by considering stratifications of arc-symmetric sets with
nonsingular strata, we prove the multiplicativity directly with formula
(\ref{r1}).

Assume therefore that $A,B$ are nonsingular arc-symmetric sets of dimension
less than or equal to $n$; suppose that $\dim A=n$ for instance.

In case $A$ compact, the result follows from a finite induction on the
dimension of $B$: indeed, resolving the singularities of $\overline B
^{\mathcal {Z}}$, one can assume that $B \subset \widetilde B$, with
$\widetilde B$ the nonsingular arc-symmetric closure of $B$. Then, by (\ref{r2})
$$\chi(A \times B)=\chi( A \times \widetilde B)-\chi\big(A \times (\widetilde B \setminus B)\big)$$
However $\chi( A \times \widetilde B)=\chi( A )\chi(\widetilde  B )$
for they are compact and nonsingular, and $$\chi\big(A \times (\widetilde B
\setminus B)\big)=\chi(A) \chi \big(\widetilde B \setminus B)\big)$$ as we can see
by stratifying $\widetilde B \setminus B$ with nonsingular strata and
using the inductive assumption of the finite induction because $\dim \widetilde B
\setminus B < \dim B$ by lemma \ref{dimen}.
Consequently $$\chi(A \times B)=\chi( A )\big(\chi(\widetilde  B )-\chi
(\widetilde B \setminus B)\big)=\chi( A )\chi( B ).$$
If $A$ is no longer compact, then compactify $A$ and $B$ in
$\widetilde A,\widetilde B$, and assume that $\widetilde A,\widetilde
B$ are nonsingular even if it means resolving singularities, as
before.

Then, by additivity, $$\chi(A \times B)=\chi(\widetilde A\times
B)-\chi \big((\widetilde A
\setminus A)\times \widetilde B\big)+\chi\big((\widetilde A
\setminus A)\times  (\widetilde B \setminus B)\big).$$
The multiplicativity of the first two terms comes from the preceding
case (in the second one, stratify the possibly singular set $\widetilde A
\setminus A$), and the multiplicativity of the third is obtained by the inductive
assumption for $\max(\dim \widetilde A \setminus A,\dim \widetilde B
\setminus B) <n$ by lemma \ref{dimen}. Therefore 
$$\chi(A \times B)=\chi(\widetilde A)\chi( B)-\chi(\widetilde A
\setminus A)\chi( \widetilde B)+\chi(\widetilde A
\setminus A)\chi (\widetilde B \setminus B)=\chi(A)\chi( B),$$
and theorem \ref{main2} is proven.
\begin{flushright} $\Box$ \end{flushright}


\section{Invariance of $\beta$ under $Nash$-isomorphisms}\label{sec3}

The definition of isomorphism between arc-symmetric sets is 
algebraic, via birational morphisms. But arc-symmetric sets are also
closely related to analytic objects. As an example, the following
proposition emphasizes the good behaviour of $\beta$ with respect
to compact algebraic varieties which are nonsingular as analytic varieties.
Recall that by $b_i (X)$ we denote the i-th Betti number of $X$ with
mod-2 coefficients, and let us put $b(X)=\sum_{i= 0}^{\dim X} b_i(X)u^i$.

\begin{prop}\label{betab} Let $X$ be a compact algebraic variety which is nonsingular as
  an analytic space. Then $\beta(X)=b(X)$.
\end{prop}
\begin{demo} One can desingularize the algebraic singularities of $X$ by
  a sequence of blowings-up with smooth centers (\cite{BM, HIRO}). At
  each step of the desingularization one has the following relation,
  where $Bl_C X$ designs the blowing-up of $X$ along the nonsingular
  subvariety $C$ with exceptionnal divisor $E$:
$$\beta(Bl_C X)-\beta(E)=\beta(X)-\beta(C)$$
because the blowing-up is birational and
$$b(Bl_C X)-b(E)=b(X)-b(C)$$
because $X$ and $C$ are smooth and the blowing-up is a degree one
morphism.
Remark that $\beta(E)=b(E)$ and $\beta(C)=b(C)$ because $E$ and $C$
are nonsingular and compact arc-symmetric sets, and the same is true for $\widetilde X$,
the desingularization of $X$. Then $\beta(X)$ and $b(X)$ can be
expressed by the same formulae in terms of $\beta$ for the former and
$b$ for the latter, where the spaces involved are nonsingular and
compact. Therefore for each one of these spaces $\beta$ and $b$
coincide, and then $\beta(X)=b(X)$.
\end{demo}

Actually we will see later that the assumption ``$X$ is an algebraic
variety'' can be replaced by the weaker ``$X$ is a semi-algebraic
set'', and so ``$X$ is a $Nash$-manifold''.

In order to relate the analytic aspect of
arc-symmetric sets and the behaviour of $\beta$, we propose the following definition of $Nash$-isomorphism
between arc-symmetric sets.

\begin{defi}\label{isoNash} Let $A,B \in \mathcal{AS}$, and assume that there exist
  compact analytic varieties $V_1,V_2$ containing $A,B$ respectively,
  and an analytic isomorphism $\phi :V_1  \longrightarrow V_2$ such
  that $\phi(A)=B$. If moreover one can choose $V_1,V_2$ to be
  semi-algebraic sets and $\phi$ a semi-algebraic map, then one say
  that  $A$ and $B$ are $Nash$-isomorphic.
\end{defi}

One wants this new definition of isomorphism respects the invariance of $\beta$ and, actually, one has the following result.

\begin{thm}\label{ASiso}  Let $A_1,A_2 \in \mathcal{AS}$ be
   $Nash$-isomorphic arc-symmetric sets. Then $\beta(A_1)=\beta(A_2)$.
\end{thm}

\begin{demo} Once more, we are going to prove the result by induction on dimension.
As a first step, let us generalize the result of proposition \ref{betab}.

\textbf{Step 1} Let $A$ be a compact arc-symmetric set which is also a
  nonsingular analytic subspace of the Zariski closure $X$ of $A$. Then $\beta(A)=b(A)$.

In order to prove this claim, one wants to apply the same method as in
the proof of proposition \ref{betab}, but if $C$ is a smooth center of blowing-up for $X$ it is not
  true in general that $C \cap A$ is nonsingular, so the equality
  $\beta(C \cap A)=b(C \cap A)$ does no longer hold. 
In order to solve this problem, consider the algebraic normalization $\widetilde
X$ of $X$. There exists $\widetilde A \subset \widetilde X$ the analytic
normalization of $A$ (\cite{loja}), which is analytically isomorphic to $A$ because
$A$ if nonsingular as an analytic space. Then $ b(A)=b(\widetilde A)$
because $b$ is invariant under homeomorphism.

Moreover
$\beta(A)=\beta(\widetilde A)$; actually the algebraic normalization is
a birational map, hence an algebraic isomorphism outside compact
subvarieties $E$ and $D$ of $\widetilde X$ and $X$ respectively, of dimension
stricly less than $\dim X=\dim A$. Thus $\beta(\widetilde A \setminus
E)=\beta(A \setminus D)$ by corollary \ref{cor-main}, and the algebraic normalization restricted to
$\widetilde A \cap E$ is an analytic isomorphism onto $A \cap D$, so
$\beta(\widetilde A \cap E)=\beta(A \cap D)$ by the inductive assumption.

Note that $\widetilde X$ is
locally analyticaly irreducible as a normal space, therefore $\widetilde A$ is a union of
connected components of $\widetilde X$. Now it is true that $C \cap \widetilde A$
is nonsingular when $C$ is nonsingular, and the method of the proof
of proposition \ref{betab} applies, therefore $\beta(\widetilde A)=b(\widetilde
A)$. It follows that $\beta(A)= \beta(\widetilde A)=b(\widetilde A)=b(A)$.

Step 1 is acheived.
\vskip 5mm
\textbf{Step 2} Let us prove the theorem in the particular case where $A_i$ is a nonsingular
arc-symmetric set for
$i=1,2$, and moreover, with assumptions of the definition of a $Nash$-isomorphism, the compact
analytic variety $V_i$ is supposed to be smooth as an analytic space, for
$i=1,2$.

\begin{itemize}
\item First we show that  $\beta (\overline {A_2}^{\mathcal{AS}})=\beta (\overline {A_1}^{\mathcal{AS}})$.

Remark that $\overline {A_2}^{\mathcal{AS}}$ is a union of
  connected components of $V_2$ by proposition
 \ref{ecl}. Thus $\overline {A_2}^{\mathcal{AS}}$ is also
 nonsingular as an analytic variety and $\beta (\overline {A_2}^{\mathcal{AS}}) =b(\overline {A_2}^{\mathcal{AS}})$ by step 1. 

Moreover  $\overline {A_2}^{\mathcal{AS}}$ is isomorphic to $\overline {A_1}^{\mathcal{AS}}$ by $\phi$:
 indeed ${\phi }^{-1}(\overline {A_2}^{\mathcal{AS}})$ is a closed arc-symmetric set because
   $\phi$ have an arc-symmetric graph and is continuous,
   and it contains $A_1$
   so $\overline {A_1}^{\mathcal{AS}}  \subset {\phi}^{-1}(\overline {A_2}^{\mathcal{AS}})$. The reverse inclusion comes from the
   fact that the image by an injective algebraic map of an
   arc-symmetric set is still an arc-symmetric set (recall that $\mathcal{AS}$
   form a constructible category, cf section 1). 
Consequently $\overline {A_1}^{\mathcal{AS}}$ is nonsingular
as an analytic variety
  because so is $\overline {A_2}^{\mathcal{AS}}$ and
    $\phi$ is an analytic isomorphism,
    hence  $\beta (\overline {A_1}^{\mathcal{AS}})=b(\overline {A_
1}^{\mathcal{AS}})$ by the first step.

Remark also that $b(\overline {A_2}^{\mathcal{AS}})=b(\overline {A_1}^{\mathcal{AS}})$ because $\phi$ is a homeomorphism between these
            two smooth compact topological varieties.

These equalities imply that
            $\beta (\overline {A_2}^{\mathcal{AS}})=\beta (\overline {A_1}^{\mathcal{AS}})$.

\item Then, remark that $\beta(\overline {A_1}^{\mathcal{AS}} \setminus
A_1)=\beta(\overline {A_2}^{\mathcal{AS}} \setminus
A_2)$. Indeed this follows from the inductive hypothesis for $\overline {A_1}^{\mathcal{AS}} \setminus
A_1$ and $\overline {A_2}^{\mathcal{AS}} \setminus
A_2$ are  $Nash$-isomorphic arc-symmetric sets of
dimension strictly less than $\dim A_2$.

\item Finally $\beta(A_1)=\beta(A_2)$. Actually
$$\beta(A_1)=\beta(\overline {A_1}^{\mathcal{AS}})-\beta(\overline {A_1}^{\mathcal{AS}} \setminus
A_1)$$
and $$\beta(A_2)=\beta(\overline {A_2}^{\mathcal{AS}})-\beta(\overline {A_2}^{\mathcal{AS}} \setminus
A_2),$$\end{itemize}
and we have proved that the second members are equal, so
$\beta(A_1)=\beta(A_2)$ by additivity of $\beta$.

\vskip 5mm
\textbf{Step 3} Reduction of the problem to Step 2.

By definition of a $Nash$-isomorphism, there exist compact analytic varieties
$V_1,V_2$ containing $A,B$, and an analytic isomorphism $\phi
 :V_1 \longrightarrow V_2$ such that $\phi(A_1)=A_2$, and moreover  $V_1$ and $V_2$ are semi-algebraic sets and $\phi$ is a semi-algebraic map.

Denote by $X_1$ and $X_2$ the Zariski closures of  $V_1$ and $V_2$.

As a first step, we are going to obtain a regular morphism rather than
a semi-algebraic map.
 Denote by $\Gamma$ the graph of $\phi$. This graph is
  semi-algebraic and analytic, thus arc-symmetric. Then the
  projection $p_i$ from $Z=\overline {\Gamma} ^{\mathcal {Z}}$ onto $X_i$, for
  $i \in \{1,2\}$, is a regular morphism whose restriction to $\Gamma$
  is an analytic isomorphism with $V_i$. Moreover the preimages by
  these restrictions of $A_1$ and $A_2$ coincide, so one can put $B=p^{-1}(A)$, where $A=A_i \subset V_i=V  \subset X_i =X$ for $i \in \{1,2\}$
and $p: Z \longrightarrow X$ denote the natural
projection. Therefore $B$ is an arc-symmetric set
 $Nash$-isomorphic to $A$, and the issue is now to prove that $\beta(B)=\beta(A)$.

In order to do this, we want to come down to nonsingular objects.

Desingularize $X$ by a sequence
of blowings-up with respect to coherent algebraic sheaves of ideals (this is
possible by \cite{BM,HIRO}). Blowing-up
$Z$ with respect to the corresponding inverse image ideal sheaves, one
has at each step a regular morphism which lifts the projection $p: Z
\longrightarrow X$ to the corresponding blowing-up by the universal property of
algebraic blowing-up.  Let $\pi_X : \widetilde X \longrightarrow
X$ denote the resolution of singularities of $X$ and $\pi_{Z} :
\widetilde {Z} \longrightarrow Z$ the corresponding
composition of blowings-up of $Z$. If $\widetilde p$ denotes the morphism
obtained between $\widetilde {Z}$ and $\widetilde X$ by universal property, one has the
following commutative diagram: 

 $$\xymatrix{ \widetilde {Z}   \ar[r]^{\widetilde p} \ar[d]_{\pi_{Z}} &  \widetilde X \ar[d]^{\pi_X}\\
     Z   \ar[r]_p  & X} $$

 Moreover $\widetilde p$ restricted to the analytic strict transform
$\widetilde {\Gamma}$ of $\Gamma$
is an analytic isomorphism onto the strict transform $\widetilde V$ of $V$ because so is
$p$ between $\Gamma$ and $V$ (here we consider the blowing-up as an
analytic one).

 $$\xymatrix{ \widetilde {\Gamma}   \ar[r]_{\sim}^{\widetilde p _{|\widetilde
       \Gamma}} \ar[d]_{{\pi_{Z}}_{|\Gamma}} &  \widetilde V \ar[d]^{{\pi_X}_{|V}}\\
     \Gamma   \ar[r]^{\sim}_{p_{|\Gamma}}  & V} $$

Now we reduce to the case where $A$ and $B$ are
nonsingular by the inductive hypothesis. Actually the singular part of $A$ and $B$ are not necessarily exchanged
by $p_{|\Gamma}$, but $\Sing(A) \cup p_{|\Gamma}^{-1}\big(\Sing(B)\big)$ and $\Sing(B)
\cup p_{|\Gamma} \big(\Sing(A)\big)$ are $Nash$-isomorphic by the restriction
of $p_{|\Gamma}$. Moreover the dimension of these arc-symmetric sets is
stricly less than $\dim A=\dim B$, so they have the
same image by $\beta$ by the inductive hypothesis. Let us denote by
$A'$ and $B'$ the complement of these sets in $A$ and $B$. So $A'$
and $B'$ are nonsingular.

As $A'$ is nonsingular, it is isomorphic, in the sense of
birational morphism, with its preimage in the desingularization
$\widetilde X$ of $X$. Consequently $B'$ is also isomorphic to its
preimage in $\widetilde {\Gamma}$ by commutativity of the first
diagram. As a consequence $\beta (A')= \beta\big(\pi ^{-1}(A')\big)$ and $\beta
(B')= \beta\big(\pi ^{-1}(B')\big)$, and we have reduced the problem
to step 2.
\end{demo}

\section{Zeta functions and blow-Nash equivalence}
\subsection{Zeta functions}\label{defzeta}
This section is directly inspired by the work of J. Denef \& F. Loeser on
their motivic zeta functions \cite{DL1}.
We first define zeta functions for an analytic germ of
functions. Then we will give a formula to compute these zeta
functions in terms of a modification.

Denote by $\mathcal L$ the space of analytic arcs at the
origin $0 \in \mathbb R ^d$:
$$\mathcal L=\mathcal L(\mathbb R ^d,0)= \{\gamma : (\mathbb R,0) \longrightarrow (\mathbb R ^d,0)
|\gamma \textrm{ is analytic}\},$$
and by $\mathcal L_n$ the space of truncated analytic arcs
$$\mathcal L_n=\mathcal L_n(\mathbb R ^d,0)= \{\gamma \in \mathcal L|
\gamma (t)=a_1t+a_2t^2+ \cdots a_nt^n,~a_i \in \mathbb R ^d\}.$$
Let $\pi_n:\mathcal L \longrightarrow \mathcal L_n$ and
$\pi_{n,i}:\mathcal L_n \longrightarrow \mathcal L_i$, with $n \geq i$, be
the truncation morphisms.

Consider $f:(\mathbb R ^d,0) \longrightarrow (\mathbb R,0)$ an
analytic function germ. We define the naive zeta function $Z_f(T)$ of $f$ as the
following element of $\mathbb Z[u,u^{-1}]((T))$: 
$$Z_f(T)= \sum _{n \geq 1}{\beta (\chi_n)u^{-nd}T^n},$$
where 
$$\chi_n =\{\gamma \in  \mathcal L_n| ord(f\circ \gamma) =n \}=\{\gamma \in  \mathcal L_n| f\circ \gamma (t)=bt^n+\cdots,
b\neq 0 \},$$
and the ``monodromic'' zeta functions, or zeta functions with signs,
as
$$Z_f^{\pm}(T)= \sum _{n \geq 1}{\beta (\chi_n^{\pm})u^{-nd}T^n},$$
where 
$$\chi_n^{\pm} =\{\gamma \in  \mathcal L_n| f\circ \gamma (t)=\pm t^n+\cdots \}.$$

Remark that $\chi_n,\chi_n^{\pm}$, for $n \geq 1$, are constructible subsets of
$\mathbb R^{nd}$, hence belong to $\mathcal{AS}$.

\begin{ex}\label{simpl} Let $f:(\mathbb R,0) \longrightarrow (\mathbb R,0)$ be defined
  by $f(x)=x^k, k \geq 1$. Then

\begin{displaymath}
\chi_n = \left\{ \begin{array}{ll}
\{ \gamma = a_mt^m+\cdots +a_nt^n|a_m \neq 0\} \simeq \mathbb R^*
\times \mathbb R^{n-m} & \textrm{if n=mk}\\
\emptyset & \textrm{otherwise}
\end{array} \right.
\end{displaymath}
Therefore $\beta (\chi_n)= (u-1)u^{n-m}$ if $n=mk$ and $0$ otherwise,
hence
$$Z_f(T)=\sum _{m \geq 1}{(u-1)u^{mk-m} \big( \frac{T}{u}
  \big)^{mk}}=(u-1)\frac{T^k}{u-T^k}.$$
To compute the zeta functions with signs, we have to consider the case
$k=2p$ and $k=2p+1$. If $k=2p$, then $\chi_n^-=\emptyset$ and

\begin{displaymath}
\chi_n^+ = \left\{ \begin{array}{ll}
\{ \gamma = \pm t^m+\cdots +a_nt^n|a_m \neq 0\} \simeq \{ \pm1\}
\times \mathbb R^{n-m} & \textrm{if n=mk}\\
\emptyset & \textrm{otherwise,}
\end{array} \right.
\end{displaymath}
so
$$Z_f^+(T)=\sum _{m \geq 1}{2u^{mk-m} \big( \frac{T}{u}
  \big)^{mk}}=2\frac{T^k}{u-T^k}.$$
If $k=2p+1$, then

\begin{displaymath}
\chi_n^{\pm} = \left\{ \begin{array}{ll}
\{ \gamma = \pm t^m+\cdots +a_nt^n|a_m \neq 0\} \simeq \{ \pm1\}
\times \mathbb R^{n-m} & \textrm{if n=mk}\\
\emptyset & \textrm{otherwise,}
\end{array} \right.
\end{displaymath}
and
$$Z_f^+(T)=Z_f^-(T)=\sum _{m \geq 1}{u^{mk-m} \big( \frac{T}{u}
  \big)^{mk}}=\frac{T^k}{u-T^k}.$$
\end{ex}

It may be convenient to express zeta functions in terms of a
modification of $f$, that is a proper birational map which is an isomorphism
over the complement of the zero locus of $f$. Actually
there exists a formula, called Denef \& Loeser formula, which enables
to do this.
In the naive case, the Denef \& Loeser formula is given by the
following proposition.
\begin{prop}\label{DLnaive} Let $\sigma:\big(M,\sigma^{-1}(0)\big)  \longrightarrow (\mathbb R ^d,0)$ be a
modification of $\mathbb R^d$ such that $f \circ \sigma$ and the
jacobian determinant $\jac \sigma$ are normal crossings simultaneously,
and assume moreover that $\sigma$ is an isomorphism over the
complement of the zero locus of $f$.

Let $(f \circ \sigma)^{-1}(0)= \cup_{j \in J}E_j$ be the decomposition
into irreducible components of $(f \circ \sigma)^{-1}(0)$, and assume
that $ \sigma^{-1}(0)=\cup_{k \in K}E_k$ for some $K \subset J$.

Put $N_i=\mult _{E_i}f \circ \sigma$ and $\nu _i=1+\mult _{E_i} \jac
\sigma$, and for $I \subset J$ denote by $E_I^0$ the set $(\cap _{i
\in I} E_i) \setminus (\cup _{j \in J \setminus I}E_j)$. Then 

$$Z_f(T)=\sum_{I\neq \emptyset} (u-1)^{|I|}\beta\big(E_I^0 \cap
\sigma^{-1}(0)\big) \prod_{i \in I}\frac{u^{-\nu_i}T^{N_i}}{1-u^{-\nu_i}T^{N_i}}.$$
\end{prop}

\begin{ex}\label{briesk2}  Let $f_k:(\mathbb R^2,0) \longrightarrow (\mathbb R,0)$ be defined
  by $f_k(x,y)=x^k+y^k$, $k\geq 2$. The blowing-up at the
  origin give a suitable modification $\sigma$ for $f$. Here $(f \circ
  \sigma)^{-1}(0)$ consists of just the exceptional divisor $\mathbb P^1$
  in case $k$ even, and furthemore, in case $k$ odd, of the strict
  transform of $f$ which is a smooth curve crossing transversally the
  exceptional divisor. Then

\begin{displaymath}
Z_{f_k} = \left\{ \begin{array}{ll}
(u^2-1)u^{-2}\frac{T^k}{1-u^{-2}T^k}   & \textrm{if k is even}\\
(u-1)\frac{u^{-2}T^k}{1-u^{-2}T^k}(u+(u-1)\frac{u^{-1}T}{1-u^{-1}T})   & \textrm{if k is odd}
\end{array} \right.
\end{displaymath}
Note in particular that for $k \neq k'$ the zeta function $Z_{f_k}$
and $Z_{f_k'}$ are different.
\end{ex}

When we are dealing with signs, one defines coverings  $\widetilde
{E_I^{0,\pm}}$ of $E_I^0$ in the following way. Let $U$ be an affine open subset of $M$
such that $f \circ \sigma=u \prod_{i\in I}y_i^{N_i}$ on $U$, with $u$
a unit. Let us put $R_{U}^{\pm}=\{ (x,t) \in (E_I^0 \cap U) \times \mathbb
R: t^m=\pm \frac{1}{u(x)}\}$, where $m=gcd(N_i)$. Then the $R_{U}^{\pm}$ glue
together along the $E_I^0 \cap U$ to give $\widetilde {E_I^{0,\pm}}$.

\begin{prop}\label{DLmono} With assumptions and notations of proposition
  \ref{DLnaive}, then

$$Z_f^{\pm}(T)=\sum_{I\neq \emptyset} (u-1)^{|I|-1}\beta\big(\widetilde{E_I^{0,\pm}} \cap
\sigma^{-1}(0)\big) \prod_{i \in I}\frac{u^{-\nu_i}T^{N_i}}{1-u^{-\nu_i}T^{N_i}}.$$
\end{prop}
We will prove these results in section 4.4.

\begin{ex}\label{exZ}
\begin{enumerate} 
\item Let $f:\mathbb R^d \longrightarrow
    \mathbb R$ be defined by $f(x)=u(x)\prod_{i=1}^k x_i^{N_i}$, with
    $N_i \in \mathbb N$. Then $f$ has already normal crossings, and 
$$Z_f(T)=(u-1)^k \prod_{i=1}^k
\frac{u^{-1}T^{N_i}}{1-u^{-1}T^{N_i}}.$$
Now, if there exists at least one $N_i$ odd, then
$Z_f^+(T)=Z_f^-(T)=\frac{1}{u-1}Z_f(T)$. On the other hand, if all $N_i$ are even, then
$Z_f^-(T)=0$ and $Z_f^+(T)=\frac{2}{u-1}Z_f(T)$ if $u$ is positive,
the converse otherwise..

\item Let $f:\mathbb R^2 \longrightarrow
    \mathbb R$ be defined by $f(x,y)=x^2+y^2$. As $f$ is positive,
    $Z_f^-(T)=0$. We obtain a
    modification by the same way as in example \ref{briesk2}, and
    $\widetilde {E_I^{0,+}}$ is here the
    boundary of a Mobius band, hence homeomorphic to $\mathbb
    P^1$. Therefore
    $Z_f^+(T)=(u+1)\frac{u^{-2}T^2}{1-u^{-2}T^2}=\frac{1}{u-1}Z_f(T)$.

\item Let $f:\mathbb R^2 \longrightarrow
    \mathbb R$ be defined by $f(x,y)=x^2+y^4$. One can solve the
    singularities of $f$ by two successive blowings-up, and one
    obtains that the exceptional divisor $E$ has 
    two irreducible components $E_1$, $E_2$ with $N_1=2, \nu_1=2,
    N_2=4, \nu _2=3$. Therefore
$$Z_f(T)=(u-1)^2\frac{u^{-2}T^2}{1-u^{-2}T^2}\frac{u^{-3}T^4}{1-u^{-3}T^4}+(u-1)u
\frac{u^{-2}T^2}{1-u^{-2}T^2}+(u-1)u \frac{u^{-3}T^4}{1-u^{-3}T^4}.$$
Moreover in this case $\widetilde {E_{\{1\}}^{0,+}}$  and
    $\widetilde {E_{\{2\}}^{0,+}}$ are homeomorphic to a circle minus
    two points, so
$$Z_f^+(T)=2(u-1)\frac{u^{-2}T^2}{1-u^{-2}T^2}\frac{u^{-3}T^4}{1-u^{-3}T^4}+(u-1)
\frac{u^{-2}T^2}{1-u^{-2}T^2}+(u-1) \frac{u^{-3}T^4}{1-u^{-3}T^4}.$$
Note that in this case one has neither $Z_f(T)=(u-1)Z_f^+(T)$ nor
$Z_f(T)=\frac{u-1}{2}Z_f^+(T)$, whereas it was the case in the
preceding examples.
    
\end{enumerate}
\end{ex}

\begin{rmk} It would be convenient to have a Thom-Sebastiani formula
  to compute the zeta functions of the function $f*g$ defined by $(x,y) \longrightarrow
  f(x)+g(y)$ from the ones of $f$ and $g$, as it is the case in
  \cite{KP,DL1,Loo}. But it seems to be impossible to
  find in general such
  formulae with our zeta functions. However in the particular
  case of two positive (respectively negative) functions, one has the
  following formulae.
\end{rmk}

\begin{prop}\label{TS} Let $f: \mathbb R^{d_1} \longrightarrow  \mathbb R$ and
  $g: \mathbb R^{d_2} \longrightarrow  \mathbb R$ be two positive or
  two negative functions. Let us put $Z_f(T)=\sum _{n \geq 1}a_nT^n$
  and $Z_g(T)=\sum _{n \geq 1}b_nT^n$, and $A_n=1-\sum_{j=1}^n a_j$,
  $B_n=1-\sum_{j=1}^n b_j$.
Then the naive zeta function of $f*g:\mathbb R^{d_1+d_2}\longrightarrow
    \mathbb R$ defined by $f*g(x,y)=f(x)+g(y)$ is $Z_{f*g}(T)=\sum _{n \geq 1}c_nT^n$ with $$c_n=a_nB_n+A_nb_n+a_nb_n.$$
\end{prop}

\begin{ex} 
\begin{enumerate}
\item Let $h:\mathbb R^2 \longrightarrow
    \mathbb R$ be defined by $h(x,y)=x^2+y^2$. Recall that
    $Z_f(T)=(u^2-1)\sum_{n \geq 1} \frac{T^{2n}}{u^{2n}}$ (cf example
    \ref{briesk2}). Putting $f(x)=g(x)=x^2$, then $h=f*g$ and by example
    \ref{simpl} we get that $a_{2n}=b_{2n}=\frac{u-1}{u^n}$ and $a_{2n+1}=b_{2n+1}=0$, hence
    $A_{2n}=A_{2n+1}=\frac{1}{u^n}$. Then by proposition \ref{TS} we rederive
    $c_{2n}=\frac{u^2-1}{u^{2n}}$ and $c_{2n+1}=0$.

\item Consider $f$ as in example \ref{exZ}.3. The odd
  coefficients are zero because $f$ and $g$ are positive, and it is
  easy to verify that $$a_{2n}=\frac{u-1}{u^n}, A_{2n}=\frac{1}{u^n},$$ and
  $$b_{4n}=\frac{u-1}{u^n}, b_{4n+2}=0,
  B_{4n}=\frac{1}{u^n}=B_{4n+2}.$$ Therefore
  $$c_{4n}=\frac{u^2-1}{u^{3n}},c_{4n+2}=\frac{u-1}{u^{3n+1}},$$ which
  was not so clear on the Denef \& Loeser formula.
\end{enumerate}
\end{ex}

\begin{demo} Remark first that $u^{nd_1}A_n=\beta(\{ \gamma \in
  \mathcal L_n:ord (f\circ \gamma) >n\})$. Actually 
$$\mathcal L_n=\pi_{n,1}^{-1}(\chi_1) \sqcup \ldots
\pi_{n,n}^{-1}(\chi_n) \sqcup \{ \gamma \in
  \mathcal L_n:ord (f\circ \gamma) >n\},$$
hence by applying $\beta$ one gets
$$u^{nd_1}=\sum_{j=1}^{n}a_iu^{nd_1} +\beta(\{ \gamma \in
  \mathcal L_n:ord (f\circ \gamma) >n\})$$ and the result follows.

Now take $(\gamma_1,\gamma_2)\in \mathcal L_n(\mathbb R^{d_1})\times
\mathcal L_n(\mathbb R^{d_2})=\mathcal L_n(\mathbb R^{d_1+d_2})$; then $ord(f \circ
\gamma_1 +g \circ \gamma_2)$ is greater than $n$ if and only if $ord (f \circ
\gamma_1)$ and $ord (g \circ \gamma_2)$ are greater than $n$ because $f$
and $g$ are of the same sign. Therefore we have to distinguish the
cases $ord (f \circ \gamma_1)=n$ and $ord (g \circ \gamma_2) >n$, $ord
(f \circ \gamma_1)>n$ and $ord (g \circ \gamma_2) =n$, and finally
$ord (f \circ \gamma_1)=n$ and $ord (g \circ \gamma_2) =n$:
$$u^{n(d_1+d_2)}\beta \big(\chi_n(f*g)\big)=\beta \big(\chi_n(f) \big)u^{nd_2}B_n+u^{nd_1}A_n\beta
\big(\chi_n(g) \big)+\beta \big(\chi_n(f) \big)\beta \big(\chi_n(g) \big).$$
\end{demo}


\subsection{Blow-$Nash$ equivalence}

Let $f,g:(\mathbb R ^d,0) \longrightarrow (\mathbb R,0)$ be two germs
of real
analytic functions. They are said to be blow-$Nash$ equivalent if there
exist two modifications $\sigma_f~:~\big(M_f,\sigma_f^{-1}(0)\big)  \longrightarrow
(\mathbb R ^d,0)$ and $\sigma_g:\big(M_g,\sigma_g^{-1}(0)\big)  \longrightarrow
(\mathbb R ^d,0)$, that is proper birational maps which are isomorphisms
over the complement of the zero locus of $f$ and $g$, and a Nash-isomorphism
(~that is a semi-algebraic map which is an analytic isomorphism) $\Phi$
between $\big(M_f,\sigma_f^{-1}(0)\big)$ and $\big(M_g,\sigma_g^{-1}(0)\big)$ which
preserves the multiplicities of the jacobian determinants of $\sigma_f$ and
$\sigma_g$ along the components of the exceptional divisors, and which
induces a homeomorphism $\phi :(\mathbb R ^d,0) \longrightarrow
(\mathbb R^d,0)$ such that $f=g \circ \phi$:

$$\xymatrix{\big(M_f,\sigma_f^{-1}(0)\big) \ar[rr]^{\Phi} \ar[d]_{\sigma_f}& &\big(M_g,\sigma_g^{-1}(0)\big) \ar[d]^{\sigma_g}\\
            (\mathbb R ^d,0)       \ar[rr]^{\phi} \ar[dr]_f        & &           (\mathbb R ^d,0) \ar[dl]_g \\
         &   (\mathbb R,0) & } $$

Blow-$Nash$ equivalence is a particular case of blow-analytic
equivalence, for which we add algebraic data. A difficult issue with
this kind of equivalence relation is to find invariants, and actually
we know just the Fukui invariants \cite{IKK} and the Zeta
functions of S. Koike and A. Parusi\'nski \cite{KP}, defined with the
Euler characteristic with compact support.

 Note that the naive zeta function generalizes the Fukui
  invariants which asign to $f$ the possible orders
  of series $f\circ \gamma$ where $\gamma$ is a real analytic arc at the
  origin. Indeed $\beta\big(\chi_n(f)\big)\neq 0$ as soon
  as $\chi_n(f) \neq \varnothing$ (recall remark \ref{degdim}), and thus the
  Fukui invariants are the nonzero exponents of $Z_f(T)$.

The following result is the main one of this section. Its proof is
directly inspired from \cite{KP}, theorem 4.5, and is a direct
consequence of lemma \ref{preDL} (see section \ref{motiv} below).

\begin{thm} The zeta functions $Z_f(T), Z_f^{\pm}(T)$ of a germ of analytic
  functions are invariants of blow-$Nash$ equivalence.  
\end{thm}

\begin{demo} Let $f,g:(\mathbb R ^d,0) \longrightarrow (\mathbb R,0)$
  be two blow-Nash equivalent analytic function germs. 
By definition there exist modifications $\sigma_f:\big(M_f,\sigma_f^{-1}(0)\big)  \longrightarrow
(\mathbb R ^d,0)$ and $\sigma_g:\big(M_g,\sigma_g^{-1}(0)\big)  \longrightarrow
(\mathbb R ^d,0)$, and a Nash-isomorphism $\Phi:
\big(M_f,\sigma_f^{-1}(0)\big) \longrightarrow \big(M_g,\sigma_g^{-1}(0)\big)$. By
\cite{HIRO,BM} one can assume that $f \circ \sigma_f$ and the jacobian
determinant $\jac \sigma_f$ (respectively $g \circ \sigma_g$ and the jacobian
determinant $\jac \sigma_g$ ) are normal crossing simultaneously,
and thus the assumptions of lemma \ref{preDL} are verified.
Now It
suffices to prove that the Denef \& Loeser formulae for the zeta functions given by
lemma \ref{preDL} coincide. 
But $\beta$ is invariant under
Nash-isomorphism by theorem \ref{ASiso}, and moreover $\phi$ preserves

-the multiplicities of $f \circ \sigma_f$ and $g \circ
\sigma_g$ because it is an isomorphism

-the multiplicities of the jacobians of $\sigma_f$ and
$\sigma_g$ along the components of the exceptional divisors by
definition of blow-Nash equivalence. 

Therefore the  zeta functions of $f$ and $g$ coincide. 
\end{demo}


\subsection{Application to Brieskorn polynomials}
We apply our zeta functions to sketch the classification of
two variables
Brieskorn polynomials under blow-$Nash$ equivalence, and to give
examples in three variables.
Brieskorn polynomials in two or three variables are polynomials of the type
$\eps_p x^p +\eps_q y^q (+\eps_r z^r)$, $ p\leq q \leq r \in
\mathbb N, \eps_p,\eps_q,\eps_r \in \{ \pm 1\}$. 

Remark that if $p=1$ then $\eps_p x +\eps_q y^q (+\eps_r z^r)$ is
$Nash$-isomorphic to $x$. Therefore we will restrict our
attention to the case $p \geq 2$.
 
Actually the
classification under blow-analytic equivalence has been done
completely in the two variable case, and almost completely in the
three variables case in \cite{KP} using zeta functions defined with Euler
characteristic with compact support and the Fukui invariants(
see \cite{KP}, theorem 7.3; for the Fukui invariant, see
\cite{IKK}). There are just one case they can not discuss, and the
following example shows that we can discuss it under blow-$Nash$
equivalence. However this it is not sufficient to conclude for blow-analytic equivalence.

\begin{ex} Let $f_{p,k}$ be the Brieskorn polynomial defined by $f_{p,k}=\pm(x^p +y^{kp}
  +z^{kp})$, $p$ even, $ k\in \mathbb N$. We are going to prove that
  for fixed $p$ and different $k$, to such polynomials are not
  blow-$Nash$ equivalent.

In order to do this we calculate
  directly the zeta function of $f_{p,k}$.
 For $n \in \mathbb N$ we compute
$\beta(\chi _n)$.
First it is clear that $\chi _n=\emptyset $ when $n$ is not a multiple
of $p$.
If $n$ is a multiple of $p$, write $n=p(mk+r)$ with $mk+r$ the
euclidean division of $\frac{n}{p}$ by $k$. If $\gamma \in
\mathcal L_n$, put $ \gamma
=(a_1t+\cdots+a_nt^n,b_1t+\cdots+b_nt^n,c_1t+\cdots+c_nt^n)$. 

Then if $r \neq 0$, the first non zero term of $f \circ \gamma$ is
given by the first component of $\gamma$, hence
$$\chi_n =\{ \gamma |a_{mk+r}\neq 0,
a_1=\cdots=a_{mk+r-1}=b_1=\cdots=b_m=c_1=\cdots=c_m=0 \}.$$
In the case where $r=0$, the three components of $\gamma$ play a part,
and 
$$\chi_n =\{ \gamma|(a_{mk},b_m,c_m)\neq (0,0,0),a_1=\cdots=a_{mk-1}=
b_1=\cdots=b_{m-1}=c_1=\cdots=c_{m-1}=0 \}.$$
Therefore

\begin{displaymath}
\chi_n \simeq \left\{ \begin{array}{ll}
 \mathbb R^* \times \mathbb R^{(p-1)(mk+r)} \times (\mathbb
  R^{p(mk+r)-m})^2 & \textrm{if $r \neq 0$}\\
(\mathbb R^3)^* \times \mathbb R^{(p-1)mk} \times (\mathbb
  R^{pmk-m})^2 & \textrm{if $r=0$}
\end{array} \right.
\end{displaymath}
hence the coefficient of $T^n$ is 

\begin{displaymath}
\beta(\chi_n)u^{-3n} = \left\{ \begin{array}{lll}
(u-1)u^{-(mk+r)-2m} & \textrm{if $n=p(mk+r), 0<r<k$}\\
(u^3-1)u^{-mk-2m} & \textrm{if $n=pmk$}\\
0 & \textrm{otherwise.}
\end{array} \right.
\end{displaymath}
Therefore the zeta function of $f_{p,k}$ looks like 
$$Z_{f_{p,k}}=(u-1) \big(
u^{-1}T^{p}+u^{-2}T^{2p}+\cdots+u^{-(k-1)}T^{(k-1)p}\big) +(u^3-1)u^{-k-2}T^{kp}$$
$$+(u-1)\big(
u^{-(k+3)}T^{(k+1)p}+u^{-(k+4)}T^{(k+2)p}+\cdots+u^{-(2k+1)}T^{(2k-1)p}\big)
+(u^3-1)u^{-2(k-2)}T^{2kp}+\cdots$$

Now it suffices to note that, for $p$ fixed and $k < k'$ the
$pk$-coefficient of $ Z_{f_{p,k}}$ is $(u^3-1)u^{-k-2}$ whereas the
one of $ Z_{f_{p,k'}}$ is $(u-1)u^{-k}$.  
\end{ex}

\begin{rmk} The case of two variables Brieskorn polynomials have been
  dealt with in \cite{KP}, using their zeta functions and the Fukui
  invariants. Actually the only case where the equivalence class of Brieskorn polynomials of
  two variables can not be distinguished using only their zeta
  functions, and which requires the use of the Fukui invariants, is the following: $f_k (x,y)=\pm (x^k+y^k)$, $k \geq 2$
  even. Remark that we have seen in example \ref{briesk2} that for $k \neq k'$ the zeta function $Z_{f_k}$
and $Z_{f_k'}$ are different, therefore our zeta function
distinguishes this case, for blow-$Nash$ equivalence.
\end{rmk}

Actually, one can say more. Indeed, for two variables
Brieskorn polynomials the naive zeta function determines the
coefficients $p$ and $q$.

\begin{prop}\label{bri2} Let $g=\pm x^p \pm y^q$ be a two variable
  Brieskorn polynomial. Then $p$ and $q$ are uniquely determined by
  the naive zeta function $Z_g(T)=\sum_{n\geq1}g_n T^n$. More precisely
$$p=min\{ n:g_n \neq 0\}$$
and if $l=min\{ n:g_n \neq a_n u^n\}$, where $\sum_{n\geq1}a_n T^n$
denote the naive zeta function of $\pm x^p$, then 

\begin{displaymath}
q = \left\{ \begin{array}{ll}
l-1 & \textrm{if $p$ id odd, $p|l-1$ and $g_l \neq (u-1) u^k$ for some
  $k\in \mathbb N$}\\
l & \textrm{otherwise.}
\end{array} \right.
\end{displaymath}
\end{prop}

\begin{demo} The characterization of $p$ is clear. Now, if $p \nmid
  q$, then $q=l$ and $g_q=(u-1)u^{k}$ for some $k\in \mathbb N$. 

If $q=kp$ for some $k\in \mathbb N$, then
  $g_{kp}=\beta(\{ \pm a^p \pm b^{kp})\})u^{2kp-k-1}$ and thus
$$g_{kp}=a_{kp}u^{kp} \iff \beta \big(\{ \pm a^p \pm b^{kp} \neq 0 \} \big)=u(u-1)
\iff p~odd.$$
In that case $g_{kp+1}=(u-1)^2 u^{2kp-k-2}$, so $l=kp+1$.

Therefore if $p$ is even $q=l$, and if $p$ is odd either $q=l$ or
$q=l-1$. More precisely if $g_l=(u-1)u^k$ then $p \nmid q$ and $q=l$,
whereas if $g_l=(u-1)^2 u^k$, then $p|l-1$ and $q=l-1$.
\end{demo}

To find the signs in front of $x^p$ and $y^q$, the naive zeta function is not
sufficient
as shown by the following example.

\begin{figure} 
\caption{Resolution tree of $x^3\pm y^4$.}\label{f3}
$$\resizebox{10cm}{!}{\rotatebox{0}{\includegraphics{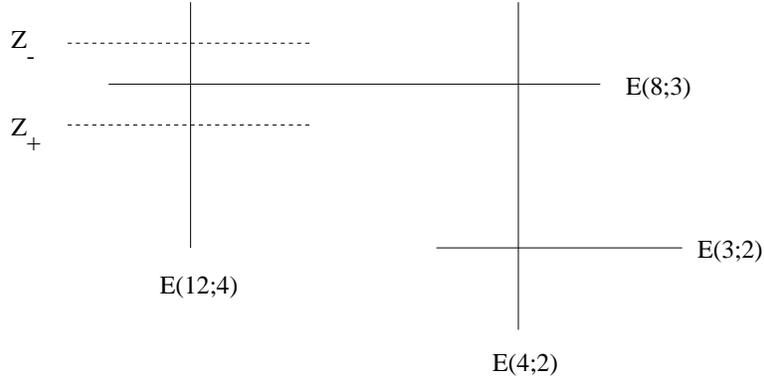}}}$$
\end{figure}

\begin{ex}  Let $f_{\pm}:(\mathbb R^2,0) \longrightarrow (\mathbb
  R,0)$ be defined by $f_{\pm}(x,y)=x^3\pm y^4$. One can solve the
  singularities of $f_{\pm}$ by a succession of four blowing-ups. The
  resolution tree of this modification $\sigma$ is drawn in figure \ref{f3},
  where $Z_{\pm}$ denotes the strict transform of $f_{\pm}$, and
  $E(N,\nu)$ denote an irreducible component of the exceptional
  divisor such that $\mult _E f_{\pm} \circ \sigma=N$ and $1+\mult _E
 \det \jac \sigma=\nu$.  
By Denef \& Loeser formula
  $Z_{f_+}=Z_{f_-}$, but $f_+$ and $f_-$ are not blow-Nash equivalent
  (they are not even blow-analytically equivalent, see \cite{KP},
  theorem 6.1 for example).   
\end{ex}

Actually the zeta functions with signs enable to discuss the signs for two
variables Brieskorn polynomials except in one
case. The followoing proposition details the possibilities.
Note
that $(x,y) \longrightarrow (\pm x,\pm y)$ give an action on the
blow-$Nash$ classes, hence when the coefficient $p$ or $q$ is odd, the
corresponding sign can not be determined. By convention, in the case
$p=q$ and the sign are opposite, we consider $x^p-y^p$ rather than $-x^p+y^p$.

\begin{prop} Let $Z_g(T)=\sum_{n\geq1}g_n T^n$ be the naive motivic
  zeta function of $\eps_p x^p+\eps_q y^q$, with $\eps_p, \eps_q \in \{\pm
  1\}$. If $p$ is even 
\begin{displaymath}
\eps_p = \left\{ \begin{array}{ll}
+1 & \textrm{if $g_p^+\neq 0$}\\
-1 & \textrm{otherwise,}
\end{array} \right.
\end{displaymath}
and if $q$ is even, but not multiple of an odd $p$, then
\begin{displaymath}
\eps_q = \left\{ \begin{array}{ll}
+1 & \textrm{if $g_p^- = 0$}\\
-1 & \textrm{otherwise.}
\end{array} \right.
\end{displaymath}
\end{prop}

\begin{rmk} In the case $p$ odd and $q=kp$ with $k$ even,
  we do not know if $x^p+y^{kp}$ and $x^p-y^{kp}$ are blow-$Nash$
  equivalent or not. They are blow-analytically equivalent by an
  argument of integration of vector field, but this method no longer holds in the
  $Nash$ setting.
\end{rmk}

\subsection{Motivic integration and the proof of Denef \& Loeser
  formula}\label{motiv}

The proof of Denef \& Loeser formula, which is a simplification of the one of \cite{DL4}, theorem 2.2.1 to
our setting, uses the theory of motivic integration on arc spaces for real algebraic
varieties (for motivic integration, we refer to \cite{DL1,DL2,DL3,DL4,Kont,Loo}). In
particular we will use the change of variables formula, which is a
fundamental result for this theory (\cite{Kont}, and \cite{DL2}, lemma 3.3). We
will recall these notions, without proofs.

Let $\sigma: \big(M, \sigma ^{-1}(0)\big) \longrightarrow (\mathbb R^d,0)$ be
a real modification.
Let 
$$\mathcal L \big(M, \sigma ^{-1}(0)\big)=\{\gamma: (\mathbb R,0)
\longrightarrow \big(M, \sigma ^{-1}(0)\big)| \gamma  \textrm{ is analytic}\}$$ be the arc space associated to $\big(M, \sigma ^{-1}(0)\big)$.
Denote by $\pi _n: \mathcal L \longrightarrow  \mathcal
  L_n$ the natural truncation morphisms, for $n \in \mathbb N$, where
  $\mathcal L$ denote eitheir $\mathcal L \big(M, \sigma ^{-1}(0)\big)$ or
  $\mathcal L (\mathbb R^d,0)$.
A subset $A \subset \mathcal L$ is called stable if $A=\pi _n^{-1}(C)$
for $C \subset \mathcal L_n$ constructible and some $n\geq 0$. Then we
can define the measure of such a stable set $A$ in $\mathbb
Z[u,u^{-1}]$ by 
$$\beta(A)=u^{-(n+1)d}\beta\big(\pi _n (A)\big)$$
for $n$ large enough ($\beta\big(\pi _n (A)\big)$ is well defined for Zariski constructible real algebraic varieties are
arc-symmetric sets). Indeed, $\beta(A)$ does not depend on $n$ because
the natural projections $\mathcal L_{n+1}
\longrightarrow  \mathcal  L_n$ are locally trivial fibrations with
fiber $\mathbb R^d$.

Let $\theta :A \longrightarrow  \mathbb Z[u,u^{-1}]$ be a map with finite image
and whose fibers are stable sets. Then the integral of $\theta$ over
$A$ with respect to $\beta$ is defined by the formula
$$\int _{A} \theta d\beta :=\sum_{c \in \mathbb Z[u,u^{-1}]}c\beta
\big(\theta^{-1}(c)\big).$$
We can state now the Kontsevich change of variables formula.

\begin{prop}\label{chgvar}(\cite{Kont,DL2}) Let $A \subset \mathcal L(\mathbb R^d,0)$ be stable, and suppose that
$\ord_t \det (\jac \sigma)$ is bounded on $\sigma ^{-1}(A)$. Then
$$\beta (A)= \int _{\sigma ^{-1}(A)} u^{-\ord_t \det (\jac \sigma)}d\beta.$$
\end{prop}

Before giving the proof of Denef \& Loeser formula, 
let us give some notations. The modification $\sigma$ induces applications $\sigma _*$
(respectively $\sigma
_{*,n}$) between $\mathcal L \big(M,\sigma^{-1}(0)\big)$ and $\mathcal
L (\mathbb R ^d,0)$ (respectively $\mathcal L_n (M,\sigma^{-1}(0)$ and $\mathcal
L_n (\mathbb R ^d,0)$). We put $\mathcal Z_n(f)=\pi_n ^{-1}(\chi_n)$ and
$\mathcal Z_n(f \circ \sigma)=\sigma_* ^{-1}\big(\mathcal Z_n(f)\big)$.
Finally for $e \geq 1$, put $\Delta _e=\{ \gamma \in \mathcal L
\big(M,\sigma^{-1}(0)\big)| \ord_t \det (\jac \sigma) \big(\gamma(t)\big)=e\}$, and
$\mathcal Z_{n,e}(f \circ \sigma)=\mathcal Z_n(f \circ \sigma) \cap
\Delta _e$.

\vskip 5mm
 Let us state a preliminary lemma.

\begin{lemma}\label{preDL} Let $\sigma:\big(M,\sigma^{-1}(0)\big)  \longrightarrow (\mathbb R ^d,0)$ be a
modification of $\mathbb R^d$ such that $f \circ \sigma$ and the
jacobian determinant $\jac \sigma$ are normal crossings simultaneously,
and assume moreover that $\sigma$ is an isomorphism over the
complement of the zero locus of $f$.

Let $(f \circ \sigma)^{-1}(0)= \cup_{j \in J}E_j$ be the decomposition
into irreducible components of $(f \circ \sigma)^{-1}(0)$, and assume
that $ \sigma^{-1}(0)=\cup_{k \in K}E_k$ for some $K \subset J$.

Put $N_i=\mult _{E_i}f \circ \sigma$ and $\nu _i=1+\mult _{E_i} \jac
\sigma$, and for $I \subset J$ denote by $E_I^0$ the set $(\cap _{i
\in I} E_i) \setminus (\cup _{j \in J \setminus I}E_j)$. Then there
exists $c\in \mathbb N$ such that
$$Z_f(T)=u^d \sum_{n \geq 1}T^n \sum_{e \leq cn}u^{-e}\sum_{I\neq
  \emptyset}\beta(\{ \gamma \in \mathcal L_n(M,E_I^0)\cap \pi_n(\Delta
_e) : f\circ \sigma
\circ \gamma (t)=bt^n+\cdots, b \neq 0\})$$
and
$$Z_f^{\pm}
(T)=u^d \sum_{n \geq 1}T^n \sum_{e \leq cn}u^{-e}\sum_{I\neq
  \emptyset}\beta(\{ \gamma \in \mathcal L_n(M,E_I^0)\cap \pi_n(\Delta
_e): f\circ \sigma
\circ \gamma (t)=\pm t^n+\cdots\}).$$

\end{lemma}

\begin{demo} Let us prove the lemma for $Z_f(T)$, the argument is the
  same for $Z_f^{\pm}$.

For $n \geq 1$, $\mathcal Z_n(f)$ is stable, so
  $\beta\big(\mathcal Z_n(f)\big)$ is defined and equals
  $u^{-(n+1)d}\beta(\chi_n)$, hence $Z_f(T)=u^d \sum_{n\geq 1}
  \beta\big(\mathcal Z_n(f)\big)T^n$.

Moreover $\mathcal Z_{n}(f \circ \sigma)$ equals the disjoint union
$\cup _{e \geq 1} \mathcal Z_{n,e}(f \circ \sigma)$. Actually this union
is finite: take $\gamma \in \mathcal Z_{n}(f \circ \sigma)$; there
exits $I \subset J$ such that $\pi _0(\gamma) \in E^0_I$. Then in a neighbourhood of $\gamma(0)$ one can choose coordinates such
that $f \circ \sigma=unit \prod_{i\in I}y_i^{N_i}$ and $\det (\jac
\sigma)=unit \prod_{i\in I}y_i^{\nu_i-1}$ where by $unit$ we denote a
non-vanishing function. Let us write
$\gamma =(\gamma_1, \ldots, \gamma_d)$, and $k_i=\ord_t \gamma_i$, for
$i=1, \ldots, d$. Then $\ord_t
f \circ \sigma \big(\gamma (t)\big)=\sum _{i=1}^d N_ik_i \geq n$ and therefore 
$$\ord_t \det (\jac \sigma) \big(\gamma(t)\big)=\sum _{i=1}^d (\nu_i-1)k_i \geq \max _i
(\frac{\nu_i-1}{N_i}) \sum _{i=1}^d N_ik_i \geq  \max
_i(\frac{\nu_i-1}{N_i})n.$$
Putting $c=\max _i (\frac{\nu_i-1}{N_i})$, we have shown thus that
$\cup _{e \geq 1} \mathcal Z_{n,e}(f \circ \sigma)=\cup _{e \leq cn}
\mathcal Z_{n,e}(f \circ \sigma)$ which is a finite union.

Now Kontsevich change of variables formula induces that
$$\beta\big(\mathcal Z_n(f)\big)=\sum_{e \leq cn} u^{-e}\beta\big(\mathcal Z_{n,e}(f \circ \sigma)\big),$$
and then
$$Z_f(T)=u^d \sum_{n\geq 1} T^n \sum_{e \leq cn} u^{-e}\beta\big(\mathcal Z_{n,e}(f \circ \sigma)\big).$$
We are going to compute $\beta\big(\mathcal Z_{n,e}(f \circ \sigma)\big)$
using the fact that $\mathcal Z_{n,e}(f \circ \sigma)$ equals the
disjoint union $$\sqcup _{I\neq
  \emptyset}  \mathcal Z_{n,e}(f \circ \sigma) \cap
\pi_0^{-1}\big(E^0_I \cap \sigma ^{-1}(0)\big),$$ thus
$$\beta\big(\mathcal Z_{n,e}(f \circ \sigma)\big)=\sum _{I\neq
  \emptyset}\beta \Big( \mathcal Z_{n,e}(f \circ \sigma) \cap
\pi_0^{-1}\big(E^0_I \cap \sigma ^{-1}(0)\big) \Big) .$$
Choose $I \neq \emptyset$. Then $\pi_n \big(\mathcal Z_{n,e}(f \circ \sigma) \cap \pi_0^{-1}(E^0_I \cap
\sigma ^{-1}(0) \big)$ is just the set
$$\{\gamma \in \mathcal L_n\big(M,\sigma
^{-1}(0)\big)|\gamma(0)\in E^0_I \cap \sigma ^{-1}(0), \ord_t f \circ
\sigma \big(\gamma(t)\big)=n, \ord_t \det (\jac \sigma) \big(\gamma(t)\big)=e \}.$$
The results follows directly from the additivity of $\beta$.
\end{demo}

The proofs of propositions \ref{DLnaive} and \ref{DLmono} just consist
of computing the value of $\beta$ on the sets which appear in the
formulae of lemma \ref{preDL}. Let us first prove proposition \ref{DLnaive}.

{\it Proof of proposition \ref{DLnaive}.}
Take $\gamma \in \pi_n \big(\mathcal Z_{n,e}(f \circ \sigma) \cap \pi_0^{-1}(E^0_I \cap
\sigma ^{-1}(0) \big)$.
On a neighbourhood of $\gamma(0)$ one can choose coordinates such
that $f \circ \sigma=unit \prod_{i\in I}y_i^{N_i}$ and moreover $\det (\jac
\sigma)=unit \prod_{i\in I}y_i^{\nu_i-1}$, hence $\pi_n(\mathcal Z_{n,e}(f \circ \sigma) \cap \pi_0^{-1}\big(E^0_I \cap
\sigma ^{-1}(0)\big)$ is isomorphic to

$$\{\gamma \in \mathcal L_n\big(M,\sigma
^{-1}(0)\big)|\gamma(0)\in E^0_I \cap \sigma ^{-1}(0), \sum_{i\in I}k_iN_i=n, \sum_{i\in I}k_i(\nu_i-1)=e\},$$
where $k_i= \ord_t \gamma_i$ for $i\in I$.
As a consequence

$$\pi_n \Big(\mathcal Z_{n,e}(f \circ \sigma) \cap \pi_0^{-1}\big(E^0_I \cap
\sigma ^{-1}(0)\big)\Big) \simeq \sqcup _{k \in A(n,e)} \big( E^0_I \cap \sigma
^{-1}(0) \big) \times (\mathbb R^*)^{|I|} \big( \prod _{i \in I} \mathbb
R^{n-k_i} \big) \times (\mathbb R^n)^{d-|I|}$$
where $A(n,e)$ is the subset of $k \in \mathbb N^d$ such that $\sum
_{i=1}^d N_ik_i=n$ and $\sum _{i=1}^d (\nu_i-1)k_i=e$.

By taking the image by $\beta$ we obtain
$$\beta\Big(\pi_n \Big(\mathcal Z_{n,e}(f \circ \sigma) \cap \pi_0^{-1}\big(E^0_I \cap
\sigma ^{-1}(0)\big)\Big)\Big)=\beta\big(E^0_I \cap \sigma
^{-1}(0)\big)(u-1)^{|I|}u^{nd-\sum _{i=1}^d k_i},$$
hence
$$Z_f(T)=\sum _{I\neq \emptyset} (u-1)^{|I|}\beta\big(E^0_I \cap \sigma
^{-1}(0)\big) \sum_{n\geq 1} \sum_{e \leq cn} \sum_{k \in A(n,e)}
u^{-e-\sum _{i=1}^d k_i}T^n.$$
Remark that $\{k \in A(n,e)|n\geq 1, e \leq cn   \}$ is in bijection
with $\mathbb N^{|I|}$,
therefore
$$ \sum_{n\geq 1} \sum_{e \leq cn} \sum_{k \in A(n,e)} u^{-e-\sum
  _{i=1}^d k_i}T^n=\sum_k \prod_{i \in
  I}(u^{-\nu_i}T^{N_i})^{k_i}=\prod_{i \in
  I}\frac{u^{-\nu_i}T^{N_i}}{1-u^{-\nu_i}T^{N_i}}.$$
Finally
$$Z_f(T)=\sum _{I\neq \emptyset} (u-1)^{|I|}\beta\big(E^0_I \cap \sigma
^{-1}(0)\big)\prod_{i \in
  I}\frac{u^{-\nu_i}T^{N_i}}{1-u^{-\nu_i}T^{N_i}},$$
which is the result.
\begin{flushright} $\Box$ \end{flushright}

The proof of proposition \ref{DLmono} is a little bit more complicated
due to the fact that we have to introduce a covering $\widetilde
{E_I^{0,\pm}}$ of $E_I^{0,\pm}$ in order to compute
$Z_f^{\pm}(T)$. Recall that if $U$ is an affine open subset of $M$
such that $f \circ \sigma=u \prod_{i\in I}y_i^{N_i}$ on $U$, with $u$
a unit, then by $R_{U}^{\pm}$ we mean $R_{U}^{\pm}=\{ (x,t) \in (E_I^0 \cap U) \times \mathbb
R: t^m=\pm \frac{1}{u(x)}\}$, where $m=gcd(N_i)$. Now $\widetilde
{E_I^{0,\pm}}$ is the gluing of the $R_{U}^{\pm}$ along the $E_I^0 \cap U$.

{\it Proof of proposition \ref{DLmono}.} Let $U$ be an affine open subset of $M$
such that $f \circ \sigma=u \prod_{i\in I}y_i^{N_i}$ on $U$, with $u$
a unit. What we have to compute is the value of $\beta$ on
$$W^{\pm}=\{ (x,y)\in (E_I^{0} \cap U)\times (\mathbb
R^*)^{|I|}:u(x)\prod_{i\in I}y_i^{N_i}=\pm1 \}.$$
Denote by $m$ the greatest common divisor of the $N_i$, $i\in I$, and
choose $n_i$, $i\in I$ such that $\sum_{i\in I} n_iN_i=m$. Assume that
$I=\{1,\ldots,s\}$.
Remark that $W^{\pm}$ is isomorphic to
$$W^{',\pm}=\{ (x,y,t)\in (E_I^{0} \cap U)\times (\mathbb
R^*)^{|I|}\times \mathbb
R^*:t^m=\frac{\pm 1}{u(x)},~\prod_{i\in I}y_i^{N_i/m}=1 \},$$
by $W^{',\pm} \longrightarrow W^{\pm}, (x,y,t) \longrightarrow
(x,t^{n_1}y_1,\ldots,t^{n_s}y_s)$. The inverse is the morphism given by $(x,y)
\longrightarrow (x, \prod_{i\in I}y_i^{N_i/m}, (\prod_{i\in
  I}y_i^{N_i/m})^{-n_1}y_1,\ldots,(\prod_{i\in
  I}y_i^{N_i/m})^{-n_s}y_s)$. 

Now it is easier to compute $\beta(W^{',\pm})$ for $W^{',\pm} \simeq R_U^{\pm} \times
(\mathbb R^*)^{|I|-1}$. This last isomorphism comes from the fact that
at least one $\frac{N_i}{m}$ is odd. Therefore
$\beta(W^{\pm})=(u-1)^{|I|-1}\beta(R_U^{\pm})$, and the same computation as in the
naive case gives the formula.
\begin{flushright} $\Box$ \end{flushright}


\end{document}